\documentclass[11pt]{amsart}
\usepackage{amssymb}
\usepackage{amsmath}
\usepackage{eurosym}
\usepackage{paralist}
\usepackage{graphics}
\usepackage{epsfig}
\usepackage{tabu}
\usepackage{graphicx}
\usepackage{epstopdf}
\usepackage[colorlinks=true]{hyperref}
\usepackage{hyperref}
\usepackage{amsmath}
\usepackage{amsfonts}
\usepackage{fancyhdr}
\usepackage{color}

\newcommand{\la}{\mathcal L_a}
\newcommand{\lat}{{\mathcal L}_{a+2}}

\newcommand{\dha}{\dot H_a^1}
\newcommand{\bd}{\mathbf{d}}
\newcommand{\dut}{\mathbf{d}(u(t))}

\newcommand{\eps}{\varepsilon}
\newcommand{\R}{\mathbb R}
\newcommand{\lsm}{\lesssim}

\newcommand{\bg}{\mathbf{g}}
\newcommand{\smt}{{\theta(t),\mu(t)}}
\newcommand{\sm}{{\theta,\mu}}
\newcommand{\lnj}{\lambda_n^j}
\newcommand{\lnk}{\lambda_n^k}
\newcommand{\xnj}{x_n^j}
\newcommand{\T}{\mathbf T}
\newtheorem{theorem}{Theorem}[section]

\newtheorem{claim}[theorem]{Claim}

\newtheorem{corollary}[theorem]{Corollary}
\newtheorem{assumption}[theorem]{Assumption}

\newtheorem{lemma}[theorem]{Lemma}

\newtheorem{proposition}[theorem]{Proposition}
\newtheorem{remark}[theorem]{Remark}

\numberwithin{equation}{section}
\DeclareMathOperator{\spa}{span}
\DeclareMathOperator{\supp}{supp}

\begin{document}
\title{Dynamics of threshold solutions for energy critical NLS with inverse square potential}

\author[K. Yang]{Kai Yang}
\address[K. Yang]{School of Mathematics, SouthEast University, Nanjing, P.R. China 211189}
\email{yangkai99sk@gmail.com, kaiyang@seu.eud.cn}

\author[C. Zeng]{Chongchun Zeng}
\address[C. Zeng]{School of Mathematics, Georgia Institute of Technology, Atlanta, GA 30332}
\email{zengch@math.gatech.edu}
\thanks{CZ was supported in part by the National Science Foundation DMS-1900083.} 

\author[X. Zhang]{Xiaoyi Zhang}
\address[X. Zhang]{Department of Mathematics, University of Iowa, Iowa City, IA 52242}
\email{xiaoyi-zhang@uiowa.edu}
\thanks{XZ was supported by Simons collaboration grant. }

\begin{abstract}

We consider the focusing energy critical NLS with inverse square potential in dimension $d= 3, 4, 5$ with the details given in $d=3$ and remarks on results in other dimensions.  Solutions on the energy surface of the ground state are characterized. We prove that solutions with kinetic energy less than that of the ground state must scatter to zero or belong to the stable/unstable manifolds of the ground state. In the latter case they converge to the ground state exponentially in the energy space as $t\to  \infty$ or $t\to -\infty$. (In 3-dim without radial assumption, this holds under the compactness assumption of non-scattering solutions on the energy surface.) When the kinetic energy is greater than that of the ground state, we show that all radial $H^1$ solutions blow up in finite time, with the only two exceptions in the case of 5-dim which belong to the stable/unstable manifold of the ground state. The proof relies on the detailed spectral analysis, local invariant manifold theory, and a global Virial analysis.

\end{abstract}
\maketitle
%\begin{center}
%{\scriptsize Keywords: dynamics, inverse square potential, NLS, energy
%critical, ground state solution.}
%\end{center}

%\address{Department of Mathematics, University of Iowa, Iowa City, IA 52242}

\section{Introduction}
Let $a\in (-\frac 14, 0)$ and $\la=-\Delta+\frac a{|x|^2}$, we consider the initial value problem 
\begin{equation*} \tag{NLS$_a$}
\begin{cases}
(i\partial_t-\la )u+|u|^4 u=0, \quad (t,x)\in \R\times\R^3, \\
u(0,x)=u_0\in \dot H^1(\R^3),
\end{cases}
\end{equation*}
for $u: \R\times\R^3\to \mathbb C$. Here the space $\dot H^1(\R^3)$ is the usual Sobolev space whose norm is given by $\|\nabla f\|_{2}$. For $a$ in the above range, the sharp Hardy's inequality implies that the bilinear form $\langle \la f, f\rangle$ is positive definite and thus defines an equivalent norm $\|\sqrt{\langle \la f, f\rangle}\|_2=\|\la^{\frac 12} f\|_2$. We use $\dha(\R^3)$ to denote the Hilbert space $\dot H^1(\R^3)$ equipped with this equivalent norm.

% induces a space $defined through the For the same reason the two norms are equivalent thus the space $\dot H^1(\R^3)$ and $\dha(\R^3)$ are isomporphism. 

The solution appearing in this paper is always a strong solution, by which we mean a function $u$ obeys the integral equation 
\begin{align*}
u(t)=e^{-it\la}u_0+i\int_0^t e^{-i(t-s)\la}|u(s)|^4u(s)ds,
\end{align*}
and lies in a certain spacetime space, for instance $u\in C_t\dot H_x^1\cap L_{t,loc}^{10}L_x^{10}$. Constructing such solution via Strichartz methodology imposes further constrains on $a$: $a>-\frac 14+\frac 1{25}$ as shown in \cite{R: Monica inverse energy}. We do not record the local theory here but would like to point out that as in the classical case, the boundedness of the spacetime norm $L_{t,x}^{10}(I\times\R^3)$ enables us to extend the solution beyond $I$ and if $I=\R$, solution scatters. Therefore we define 
%despite the fact that the full range Strichartz estimates hold for linear propagator $e^{-it\la}$. We will not record the detailed local theory. 
\[ S_I(u)=\iint_{I\times\R^3}|u(t,x)|^{10}dxdt,\]
as the scattering size of $u$. 
For a given solution $u$, we can repeatedly apply the local wellposedness to extend the solution to its maximal lifespan 
\[ (-T_*(u),T^*(u)).\]
On the interval of existence, the solution preserves its energy 
\begin{align*}
E(u(t))=\int_{\R^3}\tfrac 12|\nabla u(t,x)|^2+\tfrac{a}{2|x|^2}|u(t,x)|^2-\frac 16|u(t,x)|^6 dx.
\end{align*}
NLS${_a}$ is referred to as energy critical as the natural scaling of the equation $u(t,x)\to \lambda^{-\frac12}u(\tfrac t{\lambda^2}, \tfrac x\lambda)$ also keeps the energy invariant. 
%{\color{red}It is worthwhile to notice that the inverse square term scales the same way as Laplacian operator, this non-perturbative nature brings a lot of technical challenges to the analysis compared with the non-potential cases. }

  In the preceding work \cite{R: Harmonic inverse KMVZZ, R: Monica inverse energy}, the authors developed the fundamental analysis involving the operator $\la$ and used such to understand the scattering solutions of energy critical problem in both defocusing and focusing case.  In \cite{R: Monica inverse energy}, they proved the scattering for all finite energy solutions in the defocusing case in three dimensions and developed the crucial variational analysis of the ground state in the focusing case. Completion of the augment in multi-dimensions and focusing case  was done by the first author in \cite{R:yk_general} and \cite{R:yk_radial}. 
  
  Let us be more specific on the focusing case. In $d$ dimensions and for $a>-(\tfrac{d-2}2)^2$, the ground state soliton is the unique (up to symmetries of the equation) positive solution of static NLS$_a$: 
\begin{align}\label{groundeqn}
\mathcal L_a W=|W|^{\frac 4{d-2}}W. 
\end{align}
It was computed in \cite{R: Monica inverse energy} that  
\begin{align}\label{W}
W(x)=[d(d-2)\beta^2]^{\frac{d-2}4}\biggl(\frac{|x|^{\beta-1}}{1+|x|^{2\beta}}\biggr)^{\frac {d-2}2},\; \beta=\sqrt{1+(\tfrac 2{d-2})^2a}. 
\end{align}
Moreover, for $a\in (-(\tfrac{d-2}2)^2, 0]$, $W$ has the variational characterization which says $W$ realizes the best constant in the sharp Sobolev inequality, see for instance, \cite{R: Aubin, Lions85, R: Monica inverse energy, R: Talenti}. While for positive $a$, the problem become very tricky as the best constant can not be realized except in the radially symmetric case. We will address that case elsewhere and only focus on the case of negative $a$ in this paper.

We record the following scattering result which shows the ground state plays a role of scattering threshold. 
\begin{theorem}[\cite{R: Monica inverse energy,R:yk_general,R:yk_radial}]\label{thm:monica} 
Let $3\le d\le 6$ and $0>a> -\bigl(\tfrac{d-2}2\bigr)^2+\bigl(\tfrac{d-2}{d+2}\bigr)^2$. Let $u_0\in \dot H^1(\R^d)$ satisfy $\|u_0\|_{\dha}<\|W\|_{\dha}$ and $E(u)<E(W)$. Then there exists a unique global solution $u$ to $d$-dimensional NLS$_a$:
\[(i\partial_t-\la)u=-|u|^{\frac 4{d-2}}u, \; u(0,x)=u_0,\]
satisfying $\|u\|_{L_{t,x}^{\frac{2(d+2)}{d-2}}(\R\times\R^d)}<C(\|u_0\|_{\dha})$ in the following two scenarios: (1) $d=4,5,6$; (2) $d=3$ and $u_0$ is spherically symmetric. 
\end{theorem}
The unavailability of the result in three dimensions is ultimately due to the absence of the same scattering result for 3d quintic focusing NLS except for the spherically symmetric case. Without the radial assumption, this remains as an open problem in 3d as of now. The direct impact is the lack of compactness of non-scattering solutions on the energy surface of $E(W)$ in three dimensions. We will take the compactness as an assumption when necessary and build part of our conditional result upon it. 

Our goal in this paper is to characterize solutions on the energy surface of $E(W)$. Such problem was originated by Merle-Duyckaerts for the focusing energy critical nonlinear Schr\"odinger and wave equation in their seminal work \cite{R: D Merle, R: D Merle Wave}. We are also aware of the recent progress in \cite{SZ} on the same topic in the nonradial case. For focusing energy critical NLS, the ground state is given by the smooth bounded function 
\[
W_0(x)=\bigl(1+\tfrac{|x|^2}{d(d-2)}\bigr)^{-\frac {d-2}2},
\]
which was also proved to be the minimal energy non-scattering solution in the earlier work \cite{R: Dodson focusing NLS,R: Kenig focusing,R: Monica 5 and higher focusing}, except for $d=3$ within the class of radial data. The result in \cite{R: D Merle} demonstrated the existence of two solutions $W_0^\pm$ exponentially decaying to the ground state $W_0$ on the energy surface and classified all radial solutions as either symmetry transformations of $W_0, W_0^\pm$, scattering solutions, or blowup solutions in both time directions. While our work is largely motivated by \cite{R: D Merle}, the presence of the non-perturbative singular potential $\frac a{|x|^2}$ makes substantial differences. 
It breaks the translation symmetry of the equation and, at the same time, creates nontrivial singularity at the origin. Indeed, the fact that $\frac a{|x|^2}$ scales the same way as the Laplacian operator indicates the non-perturbative nature of this operator, making it impossible to treat the linearized problem around $W$ as a compact perturbation to any well-understood linear problem. As another example of such impact, we see the ground state $W$, which is also a stationary solution of NLS$_a$, becomes singular at the origin thus fails to belong to the full range Strichartz spaces while the free linear solutions always do \cite{R: Strichartz inverse}. As a consequence, so far even the local well-posedness of NLS$_a$ has not been established for $a$ close to $-\frac 14$. 

On the other hand, despite the disadvantage caused by the potential, the breaking of the translation symmetry 
%along with the potential 
also brings certain benefits one can take advantage of.
%, for instance, in controlling the center of concentration. I
Indeed, it has been shown in \cite{R: Monica inverse energy,R:yk_general,R:yk_radial} that the non-scattering solution on the energy surface of $E(W)$ can only concentrate around the origin instead of at any other places. Moreover, the lack of translation symmetry also indicates the manifold created by $W$ and the symmetries on the energy surface is $d$-dimension less than that in the translation invariant case. Ultimately, we are able to piece all these and the delicate spectrum analysis together to obtain the  classification of solutions on the energy surface of $E(W)$  without the radial assumption.

%Delicately taking the advantage of the state-of-art harmonic analysis related to the operator $\la$, the recent results on the spectral theory of linear Hamiltonian systems, and techniques from classical dynamical systems, we managed to obtain the classification of solutions on the energy surface of $E(W)$ even without the radial assumption. 
%

Naturally,  we need to further restrict the range of $a$ to ensure better regularity of $W$. To avoid the complexity brought up by the laborious numerology, we choose to work in dimension three even though the scattering theory in this dimension is still incomplete. Extending the 3d results to dimensions four and five is straightforward we will make a remark after each of our theorems. In the rest of higher dimensions, while most the argument can still go through, the rough nonlinearity indeed causes technical problems, for instance, in proving the Lipschitz continuity in the Strichartz spaces, a property we rely heavily on to construct the local stable/unstable manifold. Similar issue had been handled in \cite{R: LiZh NLS, R: LiZh Wave} in the case of NLS without potential. We will address the high dimension problem elsewhere.

Before stating the results, we introduce some notations. For $\theta, \mu\in \mathbb S^1\times \R^+$, we use $\T_{\theta,\mu}$ and $\bg_\sm$ to denote the symmetries transformation: 
\[\bg_\sm f(x)=e^{i\theta}\mu^{-\frac 12} f(\tfrac x\mu);\; \T_{\sm}u(t,x)=e^{i\theta}\mu^{-\frac 12}u(\tfrac t{\mu^2},\tfrac x{\mu}).\]

Our first result is the existence and uniqueness of solutions converging exponentially to $W$. 

\begin{theorem}\label{thm:main1}
Let $a\in(-\frac 14+\frac 4{25},0)$. There exist $\dot H^1(\R^3)$ solutions $W^{+}$ and $W^{-}$ to NLS$_a$ such that 
\begin{equation*}
 \lim_{t\to \infty}\|W^{\pm}(t)-W\|_{\dot H^1}\le Ce^{-ct},\;  \|W^-\|_{\dha}<\|W\|_{\dha},\; \|W^+\|_{\dha}>\|W\|_{\dha}, 
\end{equation*}
for some $C, c>0$. 
They are also unique in this class up to time translation.
Moreover, 
\begin{gather*}
W^{\pm}\in \dot H^1_{rad}(\R^3),\quad E(W^{\pm})=E(W),\\
 \int_{-\infty}^0\int_{\R^3} |W^-(t,x)|^{10}dxdt<\infty,\; W^{\pm}-W\in L^2(\R^3).
\end{gather*}
\end{theorem}
\begin{remark}
1. In dimension $d=4,5$, the same statement holds for $0>a>-(\tfrac{d-2}2)^2+(\tfrac{2(d-2)}{d+2})^2$ with the $L_{t,x}^{10}$ norm being replaced by $L_{t,x}^{\frac{2(d+2)}{d-2}}$. In particular, in dimension five where $W\in L^2(\R^5)$, $T_*(W^+)<\infty. $ See Section 7 for details. 

%{\color{blue} 2. With some additional work, it can be shown that the uniqueness holds in a larger class of solutions $\{u: \lim_{t\to \infty}\|u(t)-W\|_{\dot H^1}=0 \}$. We will address this problem in a forth coming paper. }
%
2. These solutions $W^\pm$ correspond to the two branches of the 1-dim stable manifold of $W$ in $\dot H^1 (\R^3)$, which is a smooth curve tangent to the linear stable direction at $W$. The steady state $W$ also has a 1-dim unstable manifold, given by $\overline {W^\pm}$ in this case, which satisfies the same properties in the reversed time direction. 
\end{remark}

The next result is to characterize solutions on the energy surface of $E(W)$. For the reason that was just stated, we impose the following assumption  in one part of the result. 
\begin{assumption}\label{152}
The trajectory of $\{u(t)\}$ is precompact modular scaling on $I$, i.e, there exists $\lambda(t)$ such that $\{\lambda(t)^{-\frac 12}u(t, \tfrac x{\lambda(t)}), t\in I\}$ is precompact in $\dot H^1(\R^3)$. 
\end{assumption}
We have the following
\begin{theorem}\label{thm:main2}
Let $a\in(-\frac 14+\frac 4{25},0)$. Let $u\in \dot H^1(\R^3)$ be a solution of NLS$_a$ satisfying $E(u)=E(W)$. We have

a)  If $\|u_0\|_{\dha}=\|W\|_{\dha}$, there exist $\theta,\mu$ such that $u(t,x)=\bg_\sm W$. 
%or {\color{red}$ \bg_\sm \bar W$}.

b)  If $\|u_0\|_{\dha}<\|W\|_{\dha}$, then $u$ must be a global solution. Suppose $S_\R(u)=\infty$, then $u$ conforms into one of the following two cases:\\
\indent  b.1)  $S_{[0,\infty)}(u)=\infty$. If moreover $u$ satisfies Assumption \ref{152} with $I=[0,\infty)$,  there exist $\theta,\mu,T$ such that $u(t,x)=\T_{\sm}W^-(t+T,x).$ \\
\indent  b.2)  $S_{(-\infty,0]}(u)=\infty$. If moreover $u$ satisfies Assumption \ref{152} with $I=(-\infty,0]$, then $u(t,x)=\T_{\sm}\overline {W^-(-t+T,x)}$ for some $\theta, \mu, T$. 

c) If $ \|u_0\|_{\dha}>\|W\|_{\dha}$, $u\in L^2(\R^3)$, and $u$ is radially symmetric, then $T_*(u)+T^*(u)<\infty,$  i.e. $u$ blows up both forward and backward in time.

\end{theorem}

\begin{remark}
1). Statement $b)$ in Theorem \ref{thm:main2} becomes {\bf unconditional} in four and five dimensions and in three dimensions with radial initial data. 

2.) In four dimensions, $c)$ can be stated in the same way. In five dimensions, the conclusion in $c)$ should be ``either $T_*(u)+T^*(u)<\infty$, or there exist $\theta, \mu, T$ such that $u$ equals one of the two solutions $\T_{\sm}W^+(t+T)$ and $\T_{\sm}\overline{W^+(-t+T)}$". 
\end{remark}

In the rest of the introduction we outline the main steps in the proof. 

The analysis starts with linearizing NLS$_a$ around $W$, from which we obtain a linear Hamiltonian PDE $u_t = iE''(W) u$ in the Hilbert space $\dha(\mathbb{R}^3)$ with the symplectic structure $i$ and the Hamiltonian given by the Hessian $E''(W)$ of the nonlinear energy $E(u)$. Considering $W$ is a constrained minimizer of the energy which is invariant under the phase rotation and scaling, we first prove that the quadratic form defined by $E''(W)$ has 1-dim negative direction and a 2-dim kernel based the spherical harmonics expansion and careful study on the spatial asymptotics of the resulted ODEs. Incorporating the last piece of the puzzle, i.e. the absence of the generalized kernel, we find the operator $i E''(W)$ fits right into the general framework developed in recent work \cite{lin-zeng} which immediately gives us the exponential trichotomy of $iE''(W)$. Namely, the operator $iE''(W)$ has a 1-dim stable subspace, 1-dim unstable subspace, and 1 codim-2 center subspace containing the 2-dim kernel where the linear flow has at most quadratic growth as $|t|\to \infty$. These results are summarized in Proposition \ref{prop:linzeng} in Section 3 and lays the foundation of the local nonlinear analysis of NLS$_a$. 

%The analysis starts with the linearized NLS$_a$ at $W$, which is a linear Hamiltonian PDE $u_t = iE''(W) u$ in the Hilbert space $\dha$ with the symplectic structure $i$ and the Hamiltonian given by the Hessian $E''(W)$ of the nonlinear energy $E(u)$. Since $W$ is a constrained minimizer of the energy which is invariant under the phase rotation and scaling, we first prove that the quadratic form defined by $E''(W)$ has 1-dim negative direction and a 2-dim kernel based the spherical harmonics expansion and careful study of the spatial asymptotics of the resulted ODEs. By observing that $i E''(W)$ does not have any generalized kernel, we apply the general result in \cite{lin-zeng} to obtain that $iE''(W)$ has an exponential trichotomy. Namely, it has a 1-dim stable subspace, 1-dim unstable subspace, and 1 codim-2 center subspace containing the 2-dim kernel where the linear flow has at most quadratic growth as $|t|\to \infty$. These result are summarized in Proposition \ref{prop:linzeng} in section 3 and lays the foundation of the local nonlinear analysis of NLS$_a$. 

Based on the linear analysis of $iE''(W)$, in Section 4 we establish a local coordinate near the manifold $\{ \bg_{\theta, \mu} W\}$ generated by $W$ and the symmetries. In particular, the evolution of the modulation parameter $\mu$ representing the corresponding spatial scaling would turn out to be crucial in the nonlinear analysis.

Having the exponential trichotomy decomposition from Section 3, the classical invariant manifold theory hints at the existence and uniqueness of locally invariant 1-dim stable, 1-dim unstable, and codim-2 center manifolds, see for example, \cite{CL88, JLZ18, Sc09}. To fit NLS$_a$ into the Lyapunov-Perron framework, we have to develop a Strichartz type space-time estimate for the linearized operator $iE''(W)$ with singular variable coefficients. Fortunately, treating the terms with variable coefficients as perturbations, a space-time estimate with mild temporal growth obtained by iterating a local-in-time estimate turns out to be sufficient for our construction of the local 1-dim stable manifold in Section 5. Its two branches are exactly $W^\pm$.

With the local structure being clearly established, our next step is to classify those one sided global but non-scattering solutions by proving they decay exponentially to $W$ in $\dot H^1(\mathbb{R}^3)$. Actually from the dynamical system point of view based on the saddle structure near the manifold $\{ \bg_{\theta, \mu} W\}$, such statement is rather intuitive if the solution stays in the neighborhood of this 
%manifold $\{ \bg_{\theta, \mu} W\}$ 
manifold{\footnote{In a forthcoming paper, we will show the exponential decay simply by assuming that the solution with energy $E(W)$ always stays in the neighborhood of the manifold.}, which leaves us with precluding the solution running away or traveling  into and out of small  neighborhoods. It is where the global Virial analysis comes into play. While this part of the argument is largely guided by the work in \cite{{R: D Merle}}, there are several new inputs making the proof more streamlined in the global Virial analysis.

In Section 6, we give the derivative estimate of Virial using the distance function $\dut$, which is shown to be the right quantity linking the Virial identity and the distance between $u$ and the manifold from the variational characterization of the ground state in Section 4. Solutions on the energy surface with less kinetic energy than the ground state are characterized in Section 7, where the proof of b) in Theorem \ref{thm:main2} can be found. It has been proved in the radial case and anticipated in the general case that the trajectory of such solution enjoys the precompactness after modular scaling parameter $\lambda(t)$, a property we rely heavily on in controlling the error in the Virial estimate. By properly adjusting $\lambda(t)$ (see Appendix for details), we can unify the choice of both $\lambda(t)$ and the modulation parameter $\mu(t)$ thus combine the full strength of the compactness and modulation estimates toward getting the exponential decay. The solutions on the energy surface with greater kinetic energy are considered in Section 8, where the proof of Theorem \ref{thm:main2} c) can be found. Such solutions do not have compactness, instead, we add the additional $L^2$ and radial assumption to control the error and to avoid the solution evacuating to very low frequencies. We move some of the technical estimates in the main body to the Appendix.

%With the exponential trichotomy decomposition of the linearized equation of NLS$_a$ at $W$, the classical invariant manifold theory hints at the existence and uniqueness of locally invariant 1-dim stable, 1-dim unstable, and codim-2 center manifolds, see for example, \cite{Chow-Lu-88, Schlag, Jin-Lin-Zeng}. To fit NLS$_a$ into the Lyapunov-Perron framework, we have to develop a Strichartz type space-time estimate for the linearized operator $iE''(W)$ with variable coefficients. Fortunately, treating the terms with variable coefficients as perturbations, a space-time estimate with mild temporal growth obtained by iterating a local-in-time estimate turns out to be sufficient for our construction of the local 1-dim stable manifold in Section 5. Its two branches turns out to be $W^\pm$. 
% 

\section{Preliminaries}

\textbf{Notations:} For easy reference, we include the often used notations into the following table: 

\begin{center}
\begin{tabu} to \textwidth {  X[l] |  X[l]  }
%\hline
% $A\wedge B=\min\{A,B\}$ & $A\vee B=\max\{A,B\}$  \vspace*{0.1cm}  \\
\hline
$\mathcal{L}_{a}=-\Delta +\frac{a}{|x|^{2}}$ & $\|f\|_{\dha}=\|\sqrt{\la} f\|_2$
\vspace*{0.1cm}\\
% \hline
 $\mathbf{d}(f)=|\Vert f\Vert _{\dot{H}_{a}^{1}}^{2}-\Vert W\Vert
_{\dot{H}_{a}^{1}}^{2}| $ &  $ \bg_\mu f(x)=f_{[\mu]}(x)=\mu^{-\frac 12}f(\tfrac x\mu)$
\vspace*{0.1cm}\\
%\hline
$\bg_{\sm} f(x)=f_{[\theta,\mu]}=e^{i\theta}\mu^{-\frac12}f(\tfrac x\mu)$ & $\T_{\sm}u(t,x)=e^{i\theta}\mu^{-\frac 12}u(\tfrac t{\mu^2},\tfrac x\mu)$
\vspace*{0.1cm}\\
%\hline
$\langle x\rangle=\sqrt{1+|x|^2}$ & $\beta=\sqrt{1+4a}$
\vspace*{0.1cm}\\
%\hline
$\|f\|_r=\|f\|_{L^r(\R^3)}$ & $\|f\|_{\dot H^{1,r}}=\|\sqrt{\la} f\|_r$
\vspace*{0.1cm}\\
\hline
\end{tabu}
\end{center}
%
%
%
%{\mathbf g}_\mu f=\mu ^{\frac{2-d}{2}%
%}f(\frac{x}{\mu }),\,  {\mathbf g}_{\theta,\mu}f=e^{i\theta }\mu ^{\frac{2-d}{2}%
%}f(\frac{x}{\mu })$
%\vspace*{0.1cm}\\
%\hline
%${\mathbf G}_\mu u(t,x) =\mu^{\frac{2-d}{2}}u(\frac t{\mu^2},\frac{x}{\mu })$ &
%$\langle x\rangle =\sqrt{1+|x|^{2}}$
%\vspace*{0.1cm}\\
%\hline
%   $\sigma =\tfrac{d-2}{2}-\sqrt{(\tfrac{d-2}{2})^{2}+a}$ & $\beta=\sqrt{1+\frac{4a}{(d-2)^2}}$  \vspace*{0.1cm}  \\
%  \hline
% $\|f \|_{p}=\|f \|_{L_x^p(\mathbb{R}^d)}$ & $ \|f \|_{q,r(I)}=\|f\|_{L_t^qL_x^r(I\times \mathbb{R}^d)}$  \vspace*{0.1cm}  \\
% \hline
% $\| f\| _{\dot{H}_{a}^{1}}=\Vert \sqrt{\mathcal{L}
%_{a}}f\Vert _{2} $ & $\Vert f\Vert
%_{H_{a}^{1}}=\Vert \sqrt{1+\mathcal{L}_{a}}f\Vert _{2}$ \vspace*{0.1cm}  \\
%  \hline
% $\| f\| _{\dot{H}_{a}^{1,r}}=\Vert \sqrt{\mathcal{L}
%_{a}}f\Vert _{r} $ & $\Vert f\Vert
%_{H_{a}^{1,r}}=\Vert \sqrt{1+\mathcal{L}_{a}}f\Vert _{r}$ \vspace*{0.1cm}  \\
%\hline

\textbf{Space, inner product}: Throughout this paper, we shall use $\langle \cdot, \cdot\rangle$ to denote the duality parity between a Hilbert space and its dual space. $\dha(\R^3)$ is the space of all complex functions endowed with the inner product $\Re\langle \la f,g\rangle =\Re\int\bar g \la f dx$ for any two complex functions. Occasionally, we also view $\dha(\R^3)$ as a two dimensional real valued function space and use the notation $(\dha)^2$. The same remark also applies to the Sobolev space $\dot H^1(\R^3)$.

\textbf{Variational property of the ground state $W$.}  The following lemma says $W$ is the extremizer in sharp Sobolev embedding from which one can also get the coercivity of energy. 

\begin{lemma}[\cite{R: Monica inverse energy}]\label{lem:W}
Let $a\in(-\frac 14,0)$ and $f\in \dot H^1(\R^3)$. Then 
\[\|f\|_6\le \tfrac{\|W\|_6}{\|W\|_{\dha}}\|f\|_{\dha}. \]
The equality holds if and only if $f(x)=\alpha W(\lambda x)$ for some $\alpha\in \mathbb C$ and $\lambda>0$. Moreover, if $\|f\|_{\dha}\le \|W\|_{\dha}$, then 
\begin{align}\label{12166}
\frac 13\|f\|_{\dha}^2\le E(f)\le \frac 12 \|f\|_{\dha}^2. 
\end{align}
\end{lemma}

\textbf{Strichartz estimate of $e^{-it\la}$}. We record the following linear estimate with the double endpoints estimate being given in the recent work \cite{R: Zheng Strichartz}. 

\begin{lemma}[\cite{R: Strichartz inverse,R: Zheng Strichartz}] Let $a>-\frac 14$. Let the pair of numbers $(p,q)$, $(\tilde p,\tilde q)$ satisfy 
\[\tfrac 2q+\tfrac 3r=\tfrac 2{\tilde q}+\tfrac 2{\tilde r}=\tfrac 32, \; 2\le q,\tilde q\le \infty.\]
Then the solution $u(t,x):I\times\R^3\to \mathbb C$ to the equation 
\[(i\partial_t-\la)u=f\] 
satisfy 
\begin{align*}
\|u\|_{L_t^qL_x^r(I)}\lsm \|u(t_0)\|_2+\|f\|_{L_t^{\tilde q'} L_x^{\tilde r'}(I)}
\end{align*}
for any $t_0\in I$. 
\end{lemma}
%{\color{blue}The double endpoints $(q,\tilde q)=(2,2)$ case was  confirmed in  \cite{R: Zheng Strichartz}.}

\section{Spectral analysis for the linearized operator around ground state $W$}
%Define $\mathcal L_a=-\Delta+\frac a{|x|^2}$.
%In three dimensions, let $W$ be the ground state solution of 
%\begin{align}\label{eqn}
%(iu_t-\mathcal L_a) u+|u|^4 u=0. 
%\end{align}
%i. e, 
%\begin{align}\label{groundeqn}
%\mathcal L_a W=W^{5}. 
%\end{align}
%$W$ can be calculated as 
%\begin{align}\label{W}
%W(x)=[3\beta^2]^{\frac14}\biggl(\frac{|x|^{\beta-1}}{1+|x|^{2\beta}}\biggr)^{\frac 12},
%\end{align}
%where $0>a>-\frac 14$ and 
%$$
%\beta=\sqrt{1+4a}. 
%$$

In order to study the dynamic structure of NLS$_a$ near the ground state $W$, we write the equation for $v=u-W$ in the following vector form: 
\begin{align}\label{eqnv}
\partial_t v=\mathcal L(v)+R(v). 
\end{align}
Here in the matrix form, the operator $\mathcal L$ can be written as
\begin{align*}
\mathcal L= \begin{pmatrix}
0 & \mathcal L_a-W^4\\
-\mathcal L_a+5 W^4 & 0
\end{pmatrix} 
\end{align*}
and the nonlinearity is
\begin{align*}
R(v)=i |v+W|^4(v+W) - iW^5 - 5iW^4v_1+W^4 v_2. 
\end{align*}
The linearized equation inherits the Hamiltonian structure from the nonlinear one, 
\begin{align*}
\mathcal L=JL, \ \ -i \sim J=\begin{pmatrix}0&1\\-1&0\end{pmatrix}, \mbox{ and } L=\begin{pmatrix}\mathcal L_a-5W^4&0\\
0&\mathcal L_a-W^4\end{pmatrix}
\end{align*}
where $J$ is the symplectic structure and  $L$ is the Hessian  of the energy. 
%It is worth while to notice the link between $\mathcal L$ and the diagonal matrix $L$ through 
Our first step is to understand the diagonal operator $L$ which will be further used to decode the operator $\mathcal L$ through Proposition  \ref{prop:linzeng}.

Before stating the result, we first record several facts for the operator $L:  (\dot H^1)^2 \to (\dot H^{-1})^2$, which is bounded and symmetric. 
%{\color{red}do we keep this? Throughout this section we shall use $\langle \cdot, \cdot\rangle$ to denote the duality parity between a Hilbert space and its dual space. } 
Note $W$ is the ground state solution, we have
\begin{align*}
(\la -W^4)W=0, \ (\la -5W^4)W=-4W^5<0,
\end{align*}
which implies 
\[
\langle LW, W\rangle <0. 
\]
Let $W_1$ be the generator of scaling symmetry, i.e. 
\begin{align*}
W_1=-\frac d{d\lambda}W_{[\lambda]}\biggl|_{\lambda=1}=x\cdot \nabla W+\frac 12 W.
\end{align*}
It is easy to check that
\begin{align*}
(\la-5W^4)W_1=0. 
\end{align*}

In the following lemma we will show 
%It is easy to check 
that the three directions: $W_1, iW, W$ 
%are orthogonal to each other and they 
are the only non-positive directions of $L$. 

\begin{proposition}\label{prop:L}
There exist $c,C>0$ such that the quadratic form $Q(v)=\langle Lv, v\rangle$ on $(\dot H^1)^2$ satisfies 
\[
c \|v\|_{\dot H_a^1}^2 \le Q(v) \le C\|v\|_{\dot H_a^1}^2, \quad \forall v \in X_+,
\]
where $X_+ \subset (\dot H^1)^2$ is the codim-3 closed subspace 
\[
X_+ = \{ v \in (\dot H^1)^2 \mid \langle \la W, v \rangle = \langle \la W_1, v\rangle = \langle \la (iW), v \rangle=0\}.
\]
%For the matrix operator $L$, $W$ is the only negative direction and $W_1, iW$ are the only two null directions in the sense that for any $v\in \dha$ and $v\perp Span\{W_1, iW, W\}$, the bilinear form 
%\begin{align}\label{bil}
%Q(v)=\langle Lv, v\rangle
%\end{align}
%generates a norm for such $v$:
%\begin{align}\label{normeqa}
%c\|v\|_{\dha}^2\le Q(v)\le C\|v\|_{\dha}^2.
%\end{align}

As a corollary, $L$ has one dimensional negative direction and   
\[
\ker L = \spa \{W_1, iW\}.
\]
Moreover, \[\ker (JL)^2=\ker(JL)=\ker L.\]
\end{proposition}
\begin{proof} The upper bound of $Q(v)$ follows directly from H\"older inequality and that $W\in L^6(\mathbb{R}^3)$. We will show the lower bound of $Q(v)$ by identifying the null and negative directions for each component in $L$. 

We first consider the operator $ \mathcal L_a-5W^4$ and show that there is only one negative direction in the sense that for any real scalar valued function $v\in \dot H_a^1(\mathbb{R}^3)$ and 
\begin{align}\label{ortho}
\langle \mathcal L_a v, W\rangle=0,
\end{align}
 we have 
\begin{align}\label{nonneg5}
\langle (\mathcal L_a-5W^4)v,v\rangle \ge 0. 
\end{align}
Indeed, we will see that this is an implication of the fact that $W$ is the constrained maximizer. Let 
$M=\langle \mathcal L_a W,W\rangle$ (which also equals $\int_{\mathbb{R}^3} W^6 dx$ from the ground state equation). For any $v\in \dot H^1_a(\mathbb{R}^3)$ obeying \eqref{ortho}, by taking $\mu(s)=\tfrac{M^{\frac 12}}{(M+s^2\langle \mathcal L_a v, v\rangle)^{\frac12}},$  the trajectory defined by 
\begin{align*}
l(s)=\mu(s)(W+sv)
\end{align*}
always obeys 
\begin{align*}
\langle \mathcal L_a l(s), l(s)\rangle=M. 
\end{align*}
It can be computed that 
\begin{align*}
\mu(0)=1; \mu'(0)=0; \mu''(0)=-M^{-1}\langle \mathcal L_a v, v\rangle, 
\end{align*}
and 
\begin{align*}
l(0)=W, \ l_s(0)=v, \ l_{ss}(0)=-M^{-1}\langle \mathcal L_a v, v\rangle W. 
\end{align*}
From here and noting $W$ is the constrained maximizer from Lemma \ref{lem:W}:
\begin{align}
\|W\|_6^6=\sup_{\langle \mathcal L_a w, w\rangle =M} \int_{\mathbb{R}^3} |w(x)|^6 dx, 
\end{align}
we have 
\begin{align*}
0&\ge \frac{d^2}{ds^2}\int_{\mathbb{R}^3} |l(s)|^6 dx\;\big|_{s=0}\\
&=30 \int_{\mathbb{R}^3} l(0)^4 l_s(0)^2 dx+6\int_{\mathbb{R}^3} l(0)^5 l_{ss}(0)dx\\
&=30\int_{\mathbb{R}^3} W^4 v^2 dx-6M^{-1}\langle \mathcal L_a v, v\rangle \int_{\mathbb{R}^3} W^6 dx\\
&=-6\int_{\mathbb{R}^3}(\mathcal L_a-5W^4)v\cdot v dx=-6\langle \mathcal (\la-5W^4)v,v\rangle 
\end{align*}
\eqref{nonneg5} is proved. 

Next we investigate the null direction of $L$ and it is more convenient to work in $L^2$ setting instead of $\dha$ setting. The operator $\la-5W^4$ having only one negative direction in $\dha(\mathbb{R}^3)$ implies $\la^{-\frac 12}(\la-5W^4)\la^{-\frac 12}$ has only one negative direction in $L^2(\mathbb{R}^3)$. Easily we can write 
\begin{align*}
\la^{-\frac 12}(\la-5W^4)\la^{-\frac 12}=I-5\la^{-\frac 12}W^4\la^{-\frac 12}:=I-K. 
\end{align*}
We have the following result for $K$:
\begin{claim}\label{comp}
$K:\ L^2(\R^3)\to L^2(\R^3)$ is a compact operator. 
\end{claim}

Postponing the proof for the moment, using this claim we know that $I-K$ has at most finitely many  eigenvalues in $(-\infty, \frac 12]$ which can be ordered as
$$
\lambda_1\le \lambda_2\le\cdots\le\lambda_N
$$
counting multiplicity. 

From the previous discussion and recall that
$$
(I-K)\la^{\frac 12} W_1=0,
$$
we know 
$$
\lambda_1<0, \mbox{ and }\lambda_2=0. 
$$
Our goal now is to show $\lambda_3>0$. Note as $I-K$ is symmetric we can choose eigenfunctions as the orthonormal basis of $L^2(\R^3)$ and evaluate the $L^2$ bilinear form $\langle (I-K)u, u\rangle $. Switching back to $\dha$ setting, we immediately get the desired estimate for $\la-5W^4$:
\begin{equation}\label{lb5}
\langle( \la-5W^4)u,u\rangle\ge \lambda_3 \|u\|_{\dha}^2 ,\; \forall u \perp_{\la} W, W_1.
\end{equation}

Therefore it remains to show $\lambda_3>0$ or the kernel of $I-K$ is only one-dimensional in $L^2(\R^3)$. This is equivalent to showing the kernel of $\la-5W^4$ is one dimensional in $\dot H^1(\R^3)$. 

Consider the equation 
$$
(\la-5W^4)u=0,
$$
we write $u$ in the spherical harmonic expansion:
$$
u(r,\theta)=\sum_{j=0}^\infty f_j(r)Y_j(\theta). 
$$
Here, $Y_j(\theta)$ is the $jth$ spherical harmonics and $\{ Y_j(\theta)\}_{j=0}^\infty$ form an orthonormal basis of $L^2(\mathbb S^2)$. Recall that 
\begin{gather*}
-\Delta_{\mathbb S^2} Y_j(\theta) =\mu_j Y_j(\theta), \ j=0, 1, 2,\cdots \\
0=\mu_0<\mu_1\le\mu_2\le \cdots\ \to \infty, \ Y_0=1, \mu_1=2. 
\end{gather*}
In spherical harmonic expansion, we have 
\[
(\la-5W^4)u=-\sum_{j=0}^\infty \big(( \partial_{rr} + \frac 2r \partial_r - \frac {a+\mu_j} {r^2} + 5W^4) f_j(r)\big) Y_j(\theta).
\]
Therefore we can discuss the contribution to the kernel from each spherical harmonic starting from $j=0$. 

\textbf{Case 1.} $j=0$. 

As $Y_0=1$, the kernel function in this mode must be a spherically symmetric function $u(r)$ satisfying 
$$
(\la-5W^4)u=0,
$$
which in the radial coordinate, takes the form 
\begin{align}\label{4}
u_{rr}+\frac 2r u_r+5W^4 u-\frac a{r^2} u=0.
\end{align}
Suppose $u$ is a solution independent of the known radial solution $W_1$, from Abel's theorem, we have 
\begin{align}\label{5}
u_r W_1-(W_1)_ru=\frac C{r^2}. 
\end{align}
In the small neighborhood of $r=0$, $W_1\neq 0$, we can divide both sides of \eqref{5} by $W_1^2$ and obtain,
\begin{align*}
\biggl(\frac u{W_1}\biggr)_r=\frac C{r^2W_1^2}, \ 0<r<\eps.
\end{align*}
Recalling $W_1(r)=O(r^{\frac{\beta-1}2})$ as $r\to 0^+$, integrating the above equation from $r$ to $\eps$, we have 
\begin{align*}
u(r)=O(r^{-\frac 12(1+\beta)}), \mbox{ as } r\to 0^+,
\end{align*}
which is certainly not an $\dot H^1$ function. 
Therefore $W_1$ is the unique radial kernel. 

\textbf{Case 2.} $\{j\in\mathbb N, \mu_j=2\}$. 

In this case, we assume there exists a function in the form of $G(r)Y_j(\theta)$ associated to  the $jth$ spherical harmonics in the kernel. Writing Laplacian operator in spherical coordinate,  we have
\begin{align*}
0=(\la -5W^4) (G(r)Y_j(\theta))=(\lat-5W^4)G(r) \cdot Y_j(\theta),
\end{align*}
which implies 
\begin{equation}\label{Gink}
G(r)\in \ker (\lat-5W^4). 
\end{equation}
Our first goal toward getting a contradiction is to show positivity of $G$. To this end, we take any $v\in \dot H^1(\mathbb{R}^3)$ in the spherical harmonic expansion 
$$
v:=\sum_{j=0}^\infty v_j(r)Y_j(\theta),
$$
and evaluate
\begin{align}\label{253}
\langle (\lat-5W^4)v,v\rangle =\sum_{j=0}^\infty \langle (\lat-5W^4)v_j(r),v_j(r)\rangle+\sum_{j=1}^\infty \mu_j\int_{\R^3} \frac{|v_j(x)|^2}{|x|^2} dx>0
\end{align}
As from \eqref{nonneg5}, the first summand can be estimated 
\begin{align}
\langle (\lat-5W^4)v_j(r),v_j(r)\rangle&=\langle (\lat-5W^4)v_j(r)\cdot Y_1(\theta),v_j(r)Y_1(\theta)\rangle\notag\\
&=\langle(\la-5W^4)(v_j(r)Y_1(\theta)),v_j(r)Y_1(\theta)\rangle\notag\\
&\ge 0.\label{poslat}
\end{align}
We then know that $\lat-5W^4$ is non-negative, which together with \eqref{Gink} implies that $0$ is the first eigenvalue. Hence,
$$
G(r)>0. 
$$

We now turn to looking at the equation of $G$ and $-W'$(keeping in mind that $W'<0$),
\begin{align}
&-G''-\frac 2r G'+\frac{a+2}{r^2}G-5W^4 G=0, \label{eqG}\\
&-W^{'''}-\frac 2r W''+\frac{a+2}{r^2}W'-\frac{2a}{r^3}W-5W^4W'=0.\label{eqwp}
\end{align}
Computing $[\eqref{eqG}\cdot r^2W'-\eqref{eqwp}\cdot r^2G]$, we obtain 
\begin{equation*}
r^2W'''G+2rW''G-r^2W'G''-2rW'G'+\frac{2a}r WG=0,
\end{equation*}
which can be further written into 
\begin{equation}\label{eqwG}
\frac d{dr}\bigl[r^2(W''G-W'G')\bigr]+\frac{2a}r WG=0. 
\end{equation}
Recall the asymptotics of $W$ and $G$ from \eqref{W} and Lemma \ref{lem:asymG} in Appendix:
\begin{align}\label{asym}
\begin{cases}
\mbox{As } r\to 0^+, \ G(r)=O(r^{-\frac 12+\frac 12\sqrt{9+4a}}), \ -W'=O(r^{-\frac 32+\frac \beta 2})\\
\mbox{As } r\to \infty, \ G(r)=O(r^{-\frac 12-\frac 12\sqrt{9+4a}}), \ -W'=O(r^{-\frac 32-\frac \beta 2}),
\end{cases}
\end{align}
we have
$$
-W'>G \mbox{ as } r\to 0^+\,  , \ -W'>G \mbox{ as } r\to \infty. 
$$
Let 
\[
r_0 = \sup \{ r>0\mid -W' > G \text{ on } (0, r)\}.
\]
Possibly by replacing $G$ by $CG$ for some $C>0$ sufficiently large, it holds for some $r_0 \in (0, \infty)$. 
We have  
%$r_0$ be the first point at which $-W'$ and $G$ intersect, we then have 
\begin{align*}
(W'+G)(r_0)=0, \ (W'+G)(r)<0, \ \forall r\in (0,r_0). 
\end{align*}
Hence $(W''+G')(r_0)\ge 0$ and thus 
\begin{align}\label{rl}
(W''G-W'G')(r_0)\ge 0. 
\end{align}

Using this and the positivity of $G$, we integrate \eqref{eqwG} over $(r_0,r)$ to obtain 
\begin{align}
(W''G-W'G')(r)>0,\  \forall r\in (r_0,\infty).
\end{align}
Dividing both sides by $G^2$, we have
$$
\frac d{dr}\biggl(\frac{W'}G(r)\biggr)>0, \ \forall r\in (r_0,\infty) 
$$
which in view of \eqref{asym}, contradicts with the asymptotics
\begin{align*}
\lim_{r\to \infty} \frac{W'}G(r)= -\infty
\end{align*}
for any $a\in (-\frac 14,0)$. Therefore there is no nontrivial kernel function of $\la-5W^4$ associated to the $jth$ spherical harmonics for all $j$ satisfying $\mu_j=2$. 

{\textbf{Case 3.}} $\{j\in\mathbb N, \mu_j>2\}.$

In this case, we take any function in the form of $G(r)Y_j(\theta)$, $G\neq 0$ and compute
\begin{align*}
\la (G(r)Y_j(\theta))=\lat G(r)\cdot Y_j(\theta)+\frac{\mu_j-2}{r^2}G(r)Y_j(\theta). 
\end{align*}
Using \eqref{poslat} we immediately get 
\begin{align*}
&\langle \la(G(r)Y_j(\theta)), G(r) Y_j(\theta)\rangle \\
=&\langle(\lat-5W^4) G(r), G(r)\rangle
+ (\mu_j-2)\int_{\R^3} \frac{|G(x)|^2}{|x|^2} dx>0. 
\end{align*}
This shows there is no kernel function of $\la-5W^4$ associated to $jth$ spherical harmonics for those $j$ such that $\mu_j>2$.

The positivity of $\lambda_3$ is finally proved, and we end the discussion on the operator $\la-5W^4$. 

\vspace{0.2cm}

Based on the results on $\la-5W^4$, we can get the result for $\la-W^4$ quickly. Let $\tilde \lambda_1 \le\tilde  \lambda_2 \le ...$ denote the eigenvalues of $\la^{-\frac 12} (\la-W^4) \la^{-\frac 12}$.  From 
$$
\langle (\la-5W^4)u,u\rangle<\langle (\la-W^4)u,u\rangle,
$$
we obtain $\lambda_j < \tilde \lambda_j$, $j=1,2, \ldots$. Therefore $\tilde \lambda_2 >\lambda_2 =0$ and $\tilde \lambda_1=0$ due to $\spa\{W\}= \ker (\la - W^4)$. 
%while the same proof as for $\la-5W^4$ yields $0$ , we know its kernel has dimension $1$. $W \in \ker (\la -W^4)$ and $W>0$, there must be no negative eigenvalue for $\la^{-\frac 12} (\la-W^4) \la^{-\frac 12}$ hence the first two eigenvalue must be $0=\lambda_1 <\ \tilde\lambda_2.$
This immediately implies 
$$
\langle (\la-W^4)u,u\rangle\ge \tilde \lambda_2\|u\|_{\dha}^2, \ \ \forall \mbox{ real }u \perp_{\la} W. 
$$
Combining the two parts together, we proved the estimate for $Q(v)$.

We turn to briefly proving the last statement regarding the generalized kernel. Suppose there exists a nontrivial $\dot H^1$ function $v\notin \ker (JL)$ such that \[(JL)^2 v=0.\] 
Then $v$ satisfies 
\[JLv=c_1W_1+c_2iW,\]
for some real number $c_1,c_2$ such that $c_1c_2\neq 0$.  
 Note $JL$ is a bounded operator from $(\dot H^1)^2$ to $(\dot H^{-1})^2$, we immediately get a contradiction since $W_1, iW\notin \dot H^{-1}(\mathbb{R}^3)$ as shown in the following. Take a sequence of $\dot H^1$ function with uniform norm:
\[\psi_N(x)=N^{-\frac 12}\psi(x/N),\ \psi(r)=\begin{cases} 1, &\frac 12< r\le 1\\0, & r\ge 2, \ r\le 1/4.\end{cases}\]
It is easy to see both $\int_{\R^3} W\psi_N(x)dx$ and $\int_{\R^3} W_1\psi_N(x)dx$ diverge as $N\to \infty$ by using the asymptotic estimate 
\[ W(r), W_1(r)=O(r^{-\frac12-\frac\beta2 }), \mbox{ as } r\to \infty.\]
Therefore there is no generalized kernel for $JL$.

Finally we complete the proof by verifying the Claim \ref{comp}. Indeed, note as 
$$
\la^{-\frac 12}:\ L^2(\R^3)\to \dot H^1(\R^3), \ \la^{-\frac 12}:\ \dot H^{-1}(\R^3)\to L^2(\R^3)
$$
are both bounded and the embedding $L^{\frac 65}(\R^3)\hookrightarrow \dot H^{-1}(\R^3)$ is continuous, it suffices to show 
$$
W^4\la^{-\frac 12}: \  L^2(\R^3)\to L^{\frac 65}(\R^3)
$$
is a compact operator. Taking a bounded sequence $f_n$ in $L^2(\R^3)$ and a sufficiently small number $\eps>0$, we estimate
\begin{align*}
\||\nabla|^{\varepsilon}(W^4\la^{-\frac 12}f_n)\|_{\frac 65}&\le \||\nabla|^{\varepsilon} W^4\la^{-\frac 12}f_n\|_{\frac 65}+
\|W^4|\nabla|^\varepsilon \la^{-\frac12} f_n\|_{\frac 65}\\
&\le \|W\|_{6}^3\||\nabla|^\eps W\|_{6}\|\la^{-\frac 12} f_n\|_6+\|W\|_{\frac{12}{2-\varepsilon}}^4\||\nabla|^\varepsilon \la^{-\frac12}f_n\|_{\frac 6{1+2\varepsilon}}\\
&\lesssim \|W\|_6^3\|\nabla W\|_{\frac 6{3-2\varepsilon}}\|f_n\|_2+\|W\|_{\frac {12}{2-\varepsilon}}^4\|f_n\|_2
\lesssim1. 
\end{align*}
And 
\begin{align*}
\|\chi_{>R}W^4\la^{-\frac 12}f_n\|_{\frac 65}\le \|\chi_{>R}W^4\|_{\frac 32}\|\la^{-\frac 12}f_n\|_6\lesssim R^{-c}
\end{align*}
for some positive number $c$. The compactness of $W^4\la^{-\frac 12}$ is proved, hence the Claim \ref{comp}. 
Proposition \ref{prop:L} is finally proved. 
%Finally, we briefly prove the last statement regarding the generalized kernel. Suppose there exists a nontrivial $\dot H^1$ function $v\notin \{\ker JL\}$ such that \[(JL)^2 v=0.\] Then up to a constant, $v$ satisfies \[JLv=W_1,\ \ \mbox{ or } JL v=iW.\] Note $JL$ is a bounded operator from $(\dot H^1)^2$ to $(\dot H^{-1})^2$, we immediately get a contradiction since $W, iW\notin \dot H^{-1}$ as shown in the following. Take a sequence of $\dot H^1$ function with uniform norm:
%\[\psi_N(x)=N^{-\frac 12}\psi(x/N),\ \psi(r)=\begin{cases} 1, &1< r\le 2\\0, & r\ge 2, \ r\le 1/4.\end{cases}\]
%It is easy to see both $\int_{\R^3} W\psi_N(x)dx$ and $\int_{\R^3} W_1\psi_N(x)dx$ diverge as $N\to \infty$ by using the asymptotic estimate 
%\[ W(r), W_1(r)=O(r^{-\frac12-\frac\beta2 }), \mbox{ as } r\to \infty.\]
%Therefore there is no generalized kernel for $JL$, Proposition \ref{prop:L} is finally proved. 
\end{proof}

In view of Proposition \ref{prop:L}, we are able to apply Theorem 2.1 in \cite{lin-zeng} to obtain the following

\begin{proposition}\label{prop:linzeng}
The flow $e^{tJL}$ is a well-defined operator and there exist closed subspaces $E^u,\ E^s$ and $E^c$ such that 

a) $\dim E^u=\dim E^s=1$.

b) $e^{tJL}(E^{u,s,c})=E^{u,s,c}$.

c) $ \langle Lu, u\rangle =0, \ \forall u\in E^{u,s},$ and $$E^c=\{u\in (\dot  H^1)^2; \langle Lu, v\rangle=0, \ \forall v\in E^u\oplus E^s\}.$$

%d) There exist $c,\lambda>0$ such that 
%$$
%JLu=\pm \lambda u, \forall u\in E^{u,s}, \mbox{ hence }\biggl |e^{tJL}|_{E^{u,s}}\biggr|\le Ce^{\mp \lambda t}, \ \forall t\ge 0. 
%$$

d) $\biggl |e^{tJL}\bigr|_{E^c}\biggr |\le C(1+|t|), \forall t\in \mathbb R $.

e)  $E^c=\ker L\oplus E^e$ and $(\dot H^1)^2=E^u\oplus E^s\oplus \ker L \oplus E^e$ and 

$$L\sim 
\begin{pmatrix} 0 & 1 & 0 & 0\\
1&0&0&0\\
0&0&0&0\\
0&0&0&L_e
\end{pmatrix}, 
JL\sim \begin{pmatrix}
\lambda&0&0&0\\
0&-\lambda&0&0\\
0&0&0&A_{0e}\\
0&0&0&A_e
\end{pmatrix}
$$
and 
$$
L_e\ge \epsilon>0, \ \langle L_e e^{tA_e}u, e^{tA_e}v\rangle=\langle L_e u, v\rangle.
$$         
\end{proposition}

\begin{remark}\label{R:V}
In the rest of the paper, we will assume $V^{\pm}$ is the eigenfunction taken from the $E^u$ and $E^s$: 
\[
JL V^{\pm}=\pm e_0V^{\pm}, \ e_0>0;  \textit{ and } \langle L V^+, V^-\rangle=1.
\]
We claim that $V^{\pm}\in L^2(\R^3)$. Indeed, writing $V^{\pm}=V_1\pm iV_2$, we have 
\begin{align*}
\begin{cases} 
(\la-W^4)V_2=e_0 V_1,\\
(\la-5W^4)V_1=-e_0 V_2,
\end{cases}
\end{align*}
which clearly implies 
\[e_0\int_{\R^3} |V^\pm|^2 dx=4\int_{\R^3} W^4 V_1V_2 dx\lsm \|W\|_6^4\|V_1\|_6\|V_2\|_6\lsm 1. \]
A more precise analysis on $V^\pm$ much as in Lemma \ref{lem:asymG} can be used to show that they decay exponentially in $|x|$ for sufficiently large $|x|$.  

\end{remark}

\section{Modulation analysis}
In this section, we perform the modulation analysis for solutions in the small neighborhood of the manifold $\{\bg_\sm W\}$. On energy surface of the ground state, the distance to this manifold is controlled by 
\begin{align*}
\bd(f)=\bigl | \|f\|_{\dha}^2-\|W\|_{\dha}^2\bigr|,
\end{align*}
as shown in the following result. The same result in the case of NLS can be found in \cite{R: Aubin, Lions85, R: Talenti}.

\begin{proposition}\label{P:vcW}
\label{P: VC of W}Assume that $f\in \dot{H}_{a}^{1}(\R^3)$ and $%
E(f)=E(W)$. Then for any $\varepsilon>0$, there exists $\delta>0$ such that
when 
$$
\bd(f)<\delta, \ \ \inf_{\theta\in\mathbb S^1, \mu>0}\|f-\mathbf g_{\theta,\mu} W\|_{\dha}<\varepsilon. 
$$
\end{proposition}

\begin{proof}
We argue by contradiction. Suppose the claim does not hold, then there must exist $\varepsilon _{0}>0$ and a sequence of $\dot H^1(\mathbb{R}^3)$ functions $\{f_n\}$ such that 
\begin{equation}
\text{ }E(f_{n})=E(W),\text{ }\mathbf{d}
(f_{n})\to 0,  \label{vcW 1}
\end{equation}%
but
\begin{equation}
\inf_{\theta \in \mathbb{S}^1,\mu >0}\Vert f_n-\mathbf{g}_{\theta ,\mu }W\Vert _{%
\dot{H}_{a}^{1}}>\varepsilon _{0}.  \label{vcW 2}
\end{equation}
Replacing $f_n$ by $f_n\cdot\frac{\|W\|_{\dha}}{\|f_n\|_{\dha}}$, we may assume
\begin{align}\label{guan}
\|f_n\|_{\dha}=\|W\|_{\dha}, \ \|f_n\|_6\to \|W\|_6, \ \inf_{\theta\in\mathbb S^1,\mu>0}\|f_n-\mathbf{g}_{\theta,\mu}W\|_{\dha}>\eps_0.
\end{align}
Applying Lemma \ref{lem:lpd} to $\{f_n\}$ we obtain 
$$
f_n=\sum_{j=1}^J\phi_n^j +r_n^J, 
$$
for each $J\in \{1, \cdots, J^*\}$ with the stated properties. In particular, 
%\begin{align}\label{limphi}
%\lim_{n\to \infty}\|\phi_n^j\|_{\dha}=\begin{cases}
%\|\phi\|_{\dot H^1}& \mbox{ as } \frac{|x_n^j|}{\lambda_n^j}\to \infty\\
%\|\phi\|_{\dha}& \mbox{ as } x_n^j=0. 
%\end{cases}
%\end{align}
from the $\dot H_a^1$ decoupling in Lemma \ref{lem:lpd} and \eqref{guan} we have
\begin{align}\label{933}
\|W\|_{\dha}^2= \lim_{n\to \infty}\biggl(\sum_{j=1}^J\|\phi_n^j\|_{\dha}^2+\|r_n^J\|_{\dha}^2\biggr)=\sum_{j=1}^J\|\phi^j\|_{X^j}^2+ \lim_{n\to \infty}\|r_n^J\|_{\dha}^2.
\end{align}
Here $\|\cdot\|_{X^j}=\|\cdot\|_{\dha}$ if $x_n^j\equiv 0$ and $\|\cdot\|_{X^j}=\|\cdot\|_{\dot H^1}$ if $\frac{|x_n^j|}{\lambda_n^j}\to \infty$.  As \eqref{933} holds for any $J$, we take a limit and get
\begin{align}\label{lbd}
%\limsup_{J\to J^*}
%\lim_{n\to \infty}\biggl (
\sum_{j=1}^{J^*} \|\phi^j\|_{X^j}^2
%+\|r_n^J\|_2\biggr)
\le \|W\|_{\dha}^2. 
\end{align}

On the other hand, using the decoupling in $L^6(\mathbb{R}^3)$, \eqref{guan} and the sharp Sobolev embedding, we have 
\begin{align*}
\|W\|_6^6=\lim_{n\to \infty}\|f_n\|_6^6=\sum_{j=1}^{J^*}\|\phi^j\|_6^6\le \sum_{j=1}^{J^*}\|\phi^j\|_{\dha}^6\cdot \frac{\|W\|_6^6}{\|W\|_{\dha}^6},
\end{align*}
which implies
\begin{align}\label{second}
\|W\|_{\dha}^6\le \sum_{j=1}^{J^*}\|\phi^j\|_{\dha}^6.
\end{align}
This together with \eqref{lbd} gives 
\begin{align*}
%\limsup_{J\to J^*}
\biggl(\sum_{j=1}^{J^*} \|\phi^j\|_{X^j}^2
%+\|r_n\|_{\dha}^2
\biggr)^3\le \sum_{j=1}^{J^*}\|\phi^j \|_{\dha}^6. 
\end{align*}
Note also for $a<0$, $\|\phi\|_{\dha}<\|\phi\|_{\dot H^1}$, this obviously implies that 
\[
J^*=1, \quad x_n^1 \equiv 0, \text{ and } \limsup_{n\to \infty} \|r_n^1\|_6 =0. 
\]
Therefore, \eqref{lbd} and \eqref{second} imply $\|\phi^1\|_{\dha}=\|W\|_{\dha}, \ \|\phi^1\|_6=\|W\|_6$. Moreover 
$$
f_n=(\lambda_n)^{-\frac 12}\phi^1\bigl(\frac x{\lambda_n}\bigr)+r_n^1, \mbox{  \it{and} } \|r_n^1\|_{\dha}\to 0 
$$
follow from \eqref{933}. Hence $\phi^1=\mathbf{g}_{\theta_0, \mu_0}W$ for some $\theta_0, \mu_0$. This contradicts to the last inequality in \eqref{guan}. 
\end{proof}

This together with implicit function theorem gives:

\begin{lemma} \label{lem:1st}
There exist $\delta_0, \eps_0>0$ such that for any $f\in \dha(\R^3)$ satisfying $E(f)=E(W)$ and  
 $\bd(f)<\delta_0$, there exists a unique pair $(\theta,\mu)\in \mathbb S^1\times\R^+$ such that
\[ 
\bg_{\theta,\mu}^{-1} f\perp\{iW, W_1\} \text{ and } \|f - \bg_{\theta,\mu} W\|_{\dha} < \eps_0. 
\]
Moreover, the decomposition 
\begin{equation}\label{E:decom-1} 
\bg_{\theta,\mu}^{-1} f=W+\alpha W+v, \ v\perp \{iW,W,W_1\},
\end{equation}
obeys 
\[|\alpha|\sim\|v\|_{\dha}\sim\|\bg_{\theta,\mu}^{-1}f-W\|_{\dha}\sim \bd(f).\]
\end{lemma}
 
\begin{proof} 
We prove this lemma in several steps. 

Step 1. We first focus on the neighborhood of $W$. 
Define two functionals $J_{0},J_{1}$: $\mathbb S^1\times 
\R^+\times \dot{%
H}_{a}^{1}(\R^3)\to \mathbb R$:
\begin{equation*} 
J_{0}(\theta ,\mu ,h)=\left\langle h,\mathbf g_{\theta ,\mu }(iW)\right\rangle _{\dot{%
H}_{a}^{1}},\text{ \ \ \ \ \ }J_{1}(\theta ,\mu ,h)=\left\langle h,\mathbf g_{\theta
,\mu }(W_{1})\right\rangle _{\dot{H}_{a}^{1}}.
\end{equation*}
It is easy to check that $J_0, J_1$ are linear in $h$ and $C^1$ in $\theta, \mu$. Moreover,
\begin{align*}
J_{0}(0,1,W)=J_{1}(0,1,W)=0,\ \ 
\frac{\partial(J_0, J_1)}{\partial(\theta,\mu)}\biggr|_{(0,1,W)}=
\begin{pmatrix}-\|W\|_{\dha}^2&0\\0&\|W_1\|_{\dha}^2\end{pmatrix}.
\end{align*}
Therefore the Implicit Function Theorem assures the existence of $r_1, r_2>0$ and a $C^1$ mapping $\gamma$:
$\dha(\R^3) \supset B_{r_1}(W)  \to B_{r_2}\big((0,1)\big)\subset \mathbb{S}^1\times \R^+$ such that for any $h\in B_{r_1}(W)$, 
\[
(J_0, J_1) (\theta, \mu, h) =0, \; (\theta, \mu) \in B_{r_2}\big((0,1)\big) \text{ if and only if } (\theta, \mu) = \gamma (h),
\]
which is also equivalent to 
%there exists a unique $\gamma(h)\in B_{r_2}((0,1))$ for which $J_0(\gamma(h), h)=J_1(\gamma(h),h)=0$. More specifically, for any $\|h-W\|_{\dha}\le r_1$, there exists a unique pair $|(\theta_1,\mu_1)-(0,1)|\le r_2$ such that 
$\bg_{\theta,\mu}^{-1}h\perp\{iW, W_1\}$. Moreover, due to this orthogonality, 
\[
\| h - \bg_{\theta,\mu} W\|_{\dha} = \inf_{(\theta', \mu') \subset B_{r_2}\big((0,1)\big)} \| h - \bg_{\theta',\mu'} W\|_{\dha}. 
\]

Step 2. 
%From Lemma 3.1, for any $f\in \dha(\R^3)$ satisfying $E(f)=E(W)$ and $\bd(f)<\delta_0$, there exists $(\theta_0,\mu_0)$ and a constant $\eps$ such that 
%\begin{align}\label{437}
%\|\bg_{\theta_0,\mu_0}^{-1}f-W\|_{\dha}\le \eps. 
%\end{align}
%Note $\eps\to 0$ as $\bd(f)\to 0$. Therefore by taking $\delta_0$ small enough we have $\eps<r_1$ then we are able to apply result from Step 1 directly to $\bg_{\theta_0,\mu_0}^{-1}f$ to obtain the existence of the unique pair
%\begin{align}\label{453}
%(\theta_1,\mu_1)\in\gamma(B_\eps(W))
%\end{align}
%such that for $(\theta,\mu)=(\theta_0+\theta_1,\mu_0\mu_1)$,
%\[\bg_{\theta,\mu}^{-1}f\perp \{iW,W_1\},\]
%which induces the orthogonal decomposition 
%\[\bg_{\theta,\mu}^{-1}f-W=\alpha W+v, \textit{ where } v\perp\{iW, W_1\}.\]
%On the other hand, using the continuity of the transformation $\bg_{\theta,\mu}$ in $\dha$ with respect to $\theta, \mu$, there exists another small constant $\eps_1$ such that for $(\theta_1,\mu_1)$ in \eqref{453}, 
% \[\|W-W_{[\theta_1,\mu_1]}\|_{\dha}\le \eps_1.\] This allows us to further estimate
%\[\|\bg_{\theta,\mu}^{-1}f-W\|_{\dha}\le \|\bg_{\theta_0,\mu_0}^{-1}f-W\|_{\dha}+\|W_{[\theta_1,\mu_1]}-W\|_{\dha} \le \eps+\eps_1. 
%\]
%In particular, this implies when $\bd(f)$ is sufficiently small, $\|\bg_{\theta,\mu}^{-1}f-W\|$ is sufficiently small, hence $|\alpha|$ and $\|v\|_{\dha}$ are also sufficiently small due to the orthogonality.   
We show the global uniqueness of the above pair $(\theta,\mu)$ for small $\eps_0>0$. Suppose the uniqueness is not true, then there exist $\{f_n\} \subset \dot H^1(\R^3)$ and $\{(\theta_n, \mu_n)\},\ \{(\tilde \theta_n, \tilde \mu_n)\} \subset \mathbb S^1 \times \R^+$ such that, for any $n$, $(\theta_n, \mu_n)\ne (\tilde \theta_n, \tilde \mu_n)$
\[
\bg_{\theta_n,\mu_n}^{-1} f_n, \ \bg_{\tilde \theta_n,\tilde \mu_n}^{-1} f_n \perp\{iW, W_1\}, \;\;  \|f_n - \bg_{\theta_n,\mu_n} W\|_{\dha}, \; \|f_n - \bg_{\tilde \theta_n, \tilde\mu_n} W\|_{\dha}  < \frac 1n. 
\]
This implies 
\[
\|\bg_{\theta_n - \tilde \theta_n,\mu_n/\tilde \mu_n} W - W\|_{\dha} < \frac 2n. 
\]
Recall 
\begin{claim}
Let $\{\theta_n,\mu_n\}$ be such that $\lim_{n\to \infty}\|\bg_{\theta_n,\mu_n}W-W\|_{\dha}= 0$. Then $\lim_{n\to \infty}(\theta_n,\mu_n)=(0,1)$. 
\end{claim}
The proof of this claim is a simple contradiction argument so we skip it. This contradicts the local uniqueness from which the global uniqueness follows. \\

Step 3. We prove the comparison with $\bd(f)$ under the assumption $E(f)= E(W)$. From Proposition \ref{P:vcW}, for any $\eps>0$, there exists $\delta>0$ such that for any $f\in \dha(\R^3)$ satisfying $E(f)=E(W)$ and $\bd(f)<\delta$, it holds 
\[
\inf_{(\theta', \mu') \in \mathbb S^1 \times \R^+} \| f - \bg_{\theta',\mu'} W\|_{\dha} < \eps. 
\]
By Step 1, such $f$ can be written in the form of \eqref{E:decom-1}. From the scaling invariance of energy, without loss of generality, we may consider $\theta=0, \mu=1$ only. By expanding the energy functional around $W$, we have
\begin{align*}
E(W)=&E(f)=E(W+\alpha W+v)\\
&=E(W)+\frac 12\langle E''(W)(\alpha W+v),\alpha W+v\rangle+O(\|\alpha W+v\|_{\dha}^3)\\
&=E(W)+\frac 12 \alpha^2 Q(W)+\frac 12 Q(v)+O(|\alpha|^3+\|v\|_{\dha}^3). 
\end{align*}
Here we have used the orthogonality to drop the cross term $\langle Lv, W\rangle$. This together with the ellipticity of 
$L$ on $\{W,iW,W_1\}^{\perp}$ from Lemma 2.1 and $Q(W)<0$ gives 
\[
\|v\|_{\dha}^2 \sim - \alpha^2Q(W)+ O(|\alpha|^3+\|v\|_{\dha}^3).
\]
As indicated from Step 1, if $\delta_0$ is sufficiently small, $|\alpha|$ and $\|v\|_{\dha}$ are sufficiently small accordingly, therefore we can view the cubic term as perturbation and obtain 
\[|\alpha|\sim \|v\|_{\dha}.\]
Finally, note also 
\begin{align}
\bd(f)&=\bigl|\|W+\alpha W+v\|_{\dha}^2-\|W\|_{\dha}^2\bigr|\notag\\
&=\bigl |(2\alpha+\alpha^2)\|W\|_{\dha}^2+\|v\|_{\dha}^2\bigr|=2\|W\|_{\dha}^2|\alpha|+O(\alpha^2+\|v\|_{\dha}^2).\label{924}
\end{align}
We conclude 
\[\bd(f)\sim |\alpha|\sim\|v\|_{\dha}\sim \|\alpha W+v\|_{\dha}=\|\bg_{\theta,\mu}^{-1}f-W\|_{\dha}.\]

\end{proof}
For the rest of this section, we assume $u(t)$ is a solution of $(NLS_a)$ on the time interval $I$ satisfying 
\[E(u)=E(W),\ \bd(u(t))<\delta_0, \ \forall t\in I.\]
From Lemma \ref{lem:1st}, there exists a unique pair $(\theta(t),\mu(t))$ for each $t\in I$ such that we can decompose 
\begin{align}\label{420}
\bg_{\theta(t),\mu(t)}u(t)=W+\alpha(t)W+\tilde u(t):=W+v(t), \textit{ and } \tilde u(t)\perp \{W, iW,W_1\} 
\end{align}
with $|\alpha|$ and $\| \tilde u(t)\|_{\dha}$ comparable to $\bd(u(t))$. Our next goal is to obtain the temporal derivative estimates on the modulation parameters $\theta(t)$ and $\mu(t)$. 
%Denote 
%\begin{align}\label{334}
%v(t)=g_{\theta(t),\mu(t)}u(t)-W=\alpha(t)W+\tilde u(t),
%\end{align}
%Lemma 3.2 immediately gives the following bound 
%\begin{align}\label{336}
%|\alpha(t)|\sim \|\tilde u(t)\|_{\dha}\sim\|v(t)\|_{\dha}\sim \bd( u(t))
%\end{align}

Before stating the result, we prepare a set of estimates which are needed in analyzing the modulation equation. This is where we have to trade the range of $a$ for a better integrability of the ground state $W$. 

\begin{lemma}\label{lem:error}
Let $a>-\frac 14+\frac 4{25}$. Then for any real function $R\in L^{\frac 65}(\R^3)$ and $v\in \dha(\R^3)$, we have the following bound
\begin{gather*}
|\langle R, W\rangle_{\dha}|\lsm \|R\|_{{\frac 65}}, \  |\langle v, W\rangle_{\dha}|+|\langle x\cdot \nabla v,W\rangle _{\dha}|\lsm \|v\|_{\dha},\\
|\langle \la v, W\rangle_{\dha}|+|\langle W^4v, W\rangle_{\dha}|\lsm \|v\|_{\dha}.
\end{gather*}
Here the implicit constants depend only on $W$. The same set of estimates also hold when $W$ is replaced by $W_1$. 
\end{lemma}

\begin{proof} It is straightforward to verify that under the constraint of $a$, $W, W_1\in \dot H^1(\R^3)\cap \dot H^{1,\frac{30}{11}}(\R^3)$, which by embedding, implies $W, W_1\in L^6(\R^3)\cap L^{30}(\R^3)$. 
%Moreover, 
%\[\|f\|_{\dot H^{1,p}}\sim \|f\|_{\dot H^{1,p}_a}, \textit{ when } 1<p\le \frac {30}{11}.\]
Based on these bounds we can estimate
\begin{gather*}
|\langle R, W\rangle_{\dha}|=|\langle R,\la W\rangle |=|\langle R, W^5\rangle|\le \|R\|_{\frac 65}\|W\|_{30}^5\lsm \|R\|_{\frac 65}, \\
|\langle v, W\rangle _{\dha}|\lsm \|v\|_{\dha}\|W\|_{\dha}\lsm \|v\|_{\dha},\\
\langle x\nabla v, W\rangle_{\dha}=\int_{\R^3}x\nabla v W^5 dx=-3\int_{\R^3} vW^5 dx-5\int_{\R^3}vW^4 x \nabla W dx,\\
|\langle x \nabla v, W\rangle _{\dha}|\lsm \|v\|_6\|W\|_6^5+\|v\|_6\|W\|_6^4\|x\nabla W\|_6\lsm \|v\|_{\dha},\\
|\langle \la v, W\rangle_{\dha}|=|\langle\la v, W^5\rangle|\lsm \|\nabla v\|_2\|\nabla W\|_{\frac{30}{11}}\|W\|_{30}^4\lsm \|v\|_{\dha},\\
|\langle W^4v, W\rangle_{\dha}|=|\langle W^4v, W^5\rangle |\lsm \|v\|_6\|W\|_6^4\|W\|_{30}^5\lsm \|v\|_{\dha}.
\end{gather*}
Finally as $\la W_1=5W^4 W_1$ and $W, W_1$ are both smooth functions with the same asymptotic behaviors as $|x|\to 0$ and $|x|\to \infty$, we have the same set of estimates when $W$ is replaced by $W_1$. The lemma is proved. 
\end{proof}

We are ready to state the following

\begin{lemma}\label{L: Modulation}The modulation parameters in the decomposition \eqref{420} obey
\begin{equation}
|\alpha (t)|\sim \left\Vert v(t)\right\Vert _{\dot{H}_{a}^{1}}\sim
\left\Vert \tilde{u}(t)\right\Vert _{\dot{H}_{a}^{1}}\sim \mathbf{d}(u(t)),
\label{modulation 1}
\end{equation}%
\begin{equation}
|\alpha ^{\prime }(t)|+|\theta ^{\prime }(t)|+\left\vert \frac{\mu ^{\prime
}(t)}{\mu (t)}\right\vert \lesssim\mu ^{2}(t)\mathbf{d}(u(t)).
\label{modulation 2}
\end{equation}
All the implicit constants are time independent. 
\end{lemma}
\begin{proof}
Estimate \eqref{modulation 1} follows directly from Lemma \ref{lem:1st} so we only focus on \eqref{modulation 2}. Recall 
\[u_{[\smt]}(t,x)=e^{i\theta(t)}\mu(t)^{-\frac 12}u(t,x/\mu(t)).\]
From the equation of $u$ and letting $y=x\mu(t)$ we deduce the equation for $u_{[\smt]}(t,y)$(for simplicity we drop the $t$ dependence in $\theta,\mu$ in subscript): 
\begin{align}\label{436}
(i\partial_t-\mu^2(t)\la)u_{[\sm]}+\theta'(t)u_{[\sm]}&+i\frac{\mu'(t)}{\mu(t)}(y\nabla_y+\frac 12)u_{[\sm]}\\
&=-\mu^2(t)|u_{[\sm]}|^4 u_{[\sm]}.\notag
\end{align}
Introducing the change of variable in time: $t\to s$ and $ds=\mu^2(t) dt$. Then in the $(s,y)$ variable, \eqref{436} becomes
\begin{align}\label{447}
(i\partial_s-\la)u_{[\sm]}+\theta_su_{[\sm]}+i\frac{\mu_s}{\mu}(y\nabla_y+\frac 12)u_{[\sm]}=-|u_{[\sm]}|^4u_{[\sm]}.
\end{align}
Inserting the orthogonal decomposition from \eqref{420} in $s,y$ variable: $u_{[\sm]}(s,y)=W(y)+v(s,y)$, we obtain the equation for $v:=v_1+iv_2$: 
\begin{align*}
\partial_s v+(-\la+W^4)&v_2+i(\la-5W^4)v_1\\
&-i\theta_s(v+W)+\frac{\mu_s}\mu W_1=-\frac{\mu_s}{\mu}(y\nabla_yv+\frac 12 v)+R(v).
\end{align*}
Here $R(v)$ is the high order error 
\[R(v)=i|W+v|^4(W+v)-iW^5-5iW^4v_1+W^4v_2 \]
and obeys the estimate 
\begin{align}\label{456}
\|R(v)\|_{\frac 65}\lsm \|v\|_{\dha}^2+\|v\|_{\dha}^5\lsm \bd(u(s))^2.
\end{align}
Finally inserting $v(s,y)=\alpha(s)W(y)+\tilde u(s,y)$, we obtain the equation for $\tilde u=\tilde u_1+\tilde u_2$:
\begin{align}\label{500}
\partial_s \tilde u+\alpha_s W-i\theta_s W+\frac{\mu_s}\mu W_1
+&(-\la+W^4)\tilde u_2+i(\la-5W^4)\tilde u_1-4i\alpha W^5\\
&=R(v)+i\theta_s v-\frac{\mu_s}\mu(y\nabla v+\frac 12 v). \label{509}
\end{align}
As $\tilde u\perp\{W,iW,W_1\}$, we can obtain the estimates of $\alpha_s,\theta_s,\mu_s/\mu$ simply by pairing the equation with these three directions in $\dha(\R^3)$. All the extra terms can be bounded by using Lemma \ref{lem:error} for both real and imaginary parts as showing below. 

First, we note \eqref{509} on the right side of the equation \eqref{500} only contribute the high order error. We have
\begin{align*}
&\qquad|\langle \eqref{509},W\rangle|+|\langle \eqref{509},iW\rangle|+|\langle \eqref{509},W_1\rangle|\\
&\lsm \bd(u(s))\bigl(\bd(u(s))+|\theta_s|+|\mu_s|/\mu\bigr):=\mathcal E(s). 
\end{align*}
Taking inner product between \eqref{500} and $W$, $iW$ and $W_1$ in $\dha(\R^3)$ respectively yields
\begin{align}
\alpha_s\|W\|_{\dha}^2&=\langle (\la-W^4)\tilde u_2, W\rangle_{\dha}+O(\mathcal E(s)). \\
\theta_s\|W\|_{\dha}^2&=\langle (\la-5W^4)\tilde u_1, W\rangle_{\dha}-\alpha\langle 4W^5,W\rangle_{\dha}+O(\mathcal E(s)). \\
\frac{\mu_s}\mu\|W_1\|_{\dha}^2&=\langle (\la-W^4)\tilde u_2, W_1\rangle_{\dha}+O(\mathcal E(s)). 
\end{align}
Applying Lemma \ref{lem:error} we are able to control all terms on the left sides and obtain 
\[|\alpha_s|+|\theta_s|+|\mu_s/\mu|\lsm \bd(u(s)).\]
Changing back to $t$ variable we proved \eqref{modulation 2}. The lemma is proved. 
\end{proof}

\section{Construction of local stable solutions} 

In this section, we show the existence and uniqueness of the solution converging exponentially to the ground state $W$. 

We start by proving several linear estimates of the flow $e^{tJL}$ in the Strichartz space. The way of doing it is to use the Strichartz estimate for $e^{it\la}$ and treat the $W$-related terms as perturbations. To this end, we define the Strichartz space over a time interval $I$: 
\[
\dot S^1(I)=L_t^\infty \dha\cap L_t^{5} \dot H^{1,\frac{30}{11}}(I\times\R^3). 
\]
The Sobolev norm $\|\cdot \|_{\dot H^{1, p}}$ will be estimated mostly by the operator $(\la)^{\frac12}$ 
%but works well in analyzing the nonlinear terms 
due to the equivalence of Sobolev norms developed earlier in \cite{R: Harmonic inverse KMVZZ}. The specific version we will be using is the following:

\begin{lemma}[\cite{R: Harmonic inverse KMVZZ}]
Let $a>-\frac 14+\frac 4{25}$. Then for any $p\in[\frac{30}{29}, \frac{30}{11}]$ and $f\in C_c^\infty(\R^3)$, we have 
\[\|\nabla f\|_p\sim \|(\la)^{\frac 12} f\|_p. \]
\end{lemma}

Our first estimate is about the homogeneous flow on the central space $E^c$ given in Proposition \ref{prop:linzeng}.

\begin{lemma}\label{Lemma3.1}
	Let $u_0\in E^c$ and $u(t,x)=e^{tJL} u_0$, then for any time $t> 0$, 
	\begin{align}
		\|u(\pm t)\|_{\dot H^1} &\lesssim \left\langle t \right\rangle \|u_0\|_{\dot H^1}, \label{H1bd}\\
		 \|u(\pm t)\|_2&\lsm \|u_0\|_2+\left\langle t\right\rangle^2\|u_0\|_{\dot H^1},\label{607}\\
		\|u\|_{\dot S^1([-T,T])} &\lesssim \left\langle T \right \rangle^2 \|u_0\|_{\dot H^1}. \label{local u}
	\end{align}
\end{lemma}

\begin{remark} 
Since $e^{\pm e_0t} V^\pm = e^{tJL} V^\pm$, as a corollary of the this lemma, we have $V^\pm \in \dot H^{1, \frac {30}{11}}(\R^3)$.

% A more precise analysis on $V^\pm$ will be carried out 
%%much as in Lemma \ref{lem:asymG}
%to show them decay exponentially in $|x|$ and thus belong to $L^2(\R^3)$. See Appendix for more details.
\end{remark}

\begin{proof} For simplicity we only focus on the estimate for positive times. 
     Recall from Proposition \ref{prop:linzeng}, $E^c=\ker L\oplus E^e$, and
	\begin{align*}
		JL|_{\ker L \oplus E^e} = 
		\begin{pmatrix}
		0 & A_{0e} \\
		0 & A_e
		\end{pmatrix},
	\end{align*}
	we have the expression of the linear flow
	\begin{align*}
		e^{tJL}|_{\ker L \oplus E^e}= 
		\begin{pmatrix}
			I & \int_0^t A_{0e}e^{\tau A_e}d\tau \\
			0 & e^{tA_e}
		\end{pmatrix}.
	\end{align*}
	Hence for any $u_0\in E^c$ and $
	u_0 = u_0^k + u_0^e, \mbox{ with }\ u_0^k \in \ker L, \ u_0^e \in E^e;
	$
	we can write
	\begin{align*}
		u(t)=e^{tJL}u_0 &= 
		\begin{pmatrix}
			I &  \int_0^t A_{0e}e^{\tau A_e}d\tau \\
			0 & e^{tA_e}
		\end{pmatrix}
		\begin{pmatrix}
			u_0^k \\
			u_0^e
		\end{pmatrix}\\
		&=
		\begin{pmatrix}
			u_0^k +\int_0^t A_{0e} e^{\tau A_e}u_0^e d\tau \\
			e^{tA_e}u_0^e
		\end{pmatrix}.
	\end{align*}
	From Proposition \ref{prop:linzeng} again, the second row is under control due to the ellipticity and the invariance of $L_e$: 
		$$
	\|e^{tA_e}u_0^e\|_{\dot H^1} \sim \|u_0^e\|_{\dot H^1} \lesssim \|u_0\|_{\dot H^1}.
	$$
	Plugging this estimate and using the boundedness of $A_{0e}$ on $\dot H^1(\mathbb{R}^3)$ we have
	\begin{align*}
		\|u(t)\|_{\dot H^1} \le \|e^{tA_e}u_0^e\|_{\dot H^1} + \|u_0^k\|_{\dot H^1} 
                  +\int_0^t\|A_{0e} e^{\tau A_e}u_0^e\|_{\dot H^1} d\tau 
		 \lesssim \langle t \rangle \|u_0\|_{\dot H^1}.
	\end{align*}
	\eqref{H1bd} is proved. 
	
	To prove the $L^2$ bound \eqref{607}, we use the equation of $u$
	\begin{align}\label{652}
	iu_t =\mathcal L_a u-  W^4(5u_1+iu_2).
	\end{align}
       Multiplying both sides by $\bar u$, taking the imaginary part and integrating over $[0,t]\times \mathbb{R}^3$ gives
       \begin{align*}
       \|u(t)\|_2^2-\|u_0\|_2^2=8 \int_0^t\int_{\R^3} W^4 u_1u_2 dx ds\lsm  t \|W\|_6^4\|u\|_{L_t^\infty \dot H^1([0,t])}^2
       \end{align*}
       which together with the $\dot H^1 $ estimate from \eqref{H1bd} yields \eqref{607}.

	We turn to the estimate \eqref{local u}. Take a small number $\eta$ we partition$[0,T]$ into
	$$
	[0, T] = \bigcup_{j=0}^{N} I_j, \mbox{ with } I_j=[j\eta, (j+1)\eta], \ j\le N-1; I_{N}=[N\eta,T].
	$$
	On each interval $I_j$, by the Strichartz estimate of $e^{it\la}$ in \cite{R: Strichartz inverse}, we obtain
%	$u(t,x)$, which in the complex form, solves 
%	$$
%	iu_t =\mathcal L_a u-  W^4(5u_1+iu_2).
%	$$
%By the Strichtartz estimate for $e^{it\la}$ proved in \cite{R: Monica inverse energy}, we obtain
	\begin{align*}
		\|u\|_{\dot S^1(I_j)}&\lesssim \|u(j\eta)\|_{\dot H^1} +
		    \|\mathcal L_a^{\frac 12} (W^4(5u_1+iu_2))\|_{L^1_tL_x^2(I_j)} \\
		 &\lesssim \langle j\eta \rangle \|u_0\|_{\dot H^1} + 
		   \sum_{i=1}^2 \|W^4\nabla u_i\|_{L_t^1L_x^2(I_j)} + \|\nabla W W^3 u\|_{L_t^1L_x^2(I_j)} \\
		 &\lsm  \langle j\eta \rangle \|u_0\|_{\dot H^1} +
		 \eta^{\frac 45}(\|W\|_{{30}}^4 \|\nabla u\|_{L_t^5L_x^{\frac {30}{11}}(I_j)} 
		 + \|W\|_{30}^3\|\nabla W\|_{{\frac {30}{11}}}\|u\|_{L_t^5L_x^{30}(I_j)}) \\
		 &\lesssim  \langle j\eta \rangle \|u_0\|_{\dot H^1} + \eta^{\frac 45}\|u\|_{\dot S^1(I_j)}.
	\end{align*}
	Taking $\eta$ sufficiently small, we obtain
	$$
	\|u\|_{\dot S^1([j\eta, (j+1)\eta))} \lsm\langle j\eta \rangle \|u_0\|_{\dot H^1}.
	$$
	Summing in $j$ we obtain \eqref{local u}.			
\end{proof}

Next we prove the estimate for the inhomogeneous term.

\begin{lemma}\label{lemma3.2}
	Let $f\in E^c$ and
	$$
	v(t,x)=\int_0^t e^{(t-s)JL}f(s)ds,\ w(t,x)=\int_t^T e^{(t-s)JL}f(s)ds,
	$$
	then
	\begin{align} 
		\|v\|_{\dot S^1([0,T])}+\|w\|_{\dot S^1([0,T])}\lesssim \langle T\rangle^2 \|f\|_{L_t^1 \dot H^1([0,T])}.\label{inhomo}
	\end{align}
\end{lemma}
\begin{proof} We only prove the estimate for $v$ as the other one is similar. 
	Again, we  partition the interval $[0,T]$ into subintervals as in Lemma \ref{Lemma3.1} and apply the Strichartz estimate on 
	$I_j=[j\eta, (j+1)\eta]$ to $v(t,x)$ which solves
	$$
	iv_t=\mathcal L_a v-(5W^4v_1+iW^4v_2)+if.
	$$ 
	We have
	\begin{align*}
		\|v\|_{\dot S^1(I_j)} &\lesssim \|v(j\eta)\|_{\dot H^1} + \|5W^4 v_1 +iW^4v_2\|_{L_t^1\dot H^1(I_j)}
                  +\|f\|_{L_t^1\dot H^1(I_j)} \\
		&\lesssim \|v(j\eta)\|_{\dot H^1} +\eta^{\frac 45}\|v\|_{\dot S^1(I_j)} +\|f\|_{L_t^1\dot H^1(I_j)}
	\end{align*}
	Taking $\eta$ small enough and using \eqref{H1bd} from Lemma \ref{Lemma3.1}, we have
	\begin{align*}
		\|v(t)\|_{\dot H^1} &\le \int_0^t\|e^{(t-s)JL}f(s)\|_{\dot H^1} ds 
		 \lesssim \int_0^t\langle t \rangle \|f(s)\|_{\dot H^1} ds 
		\le \langle t \rangle \|f\|_{L_s^1\dot H^1([0,t])}.
	\end{align*}
	From here, we continue the estimate of $v$ and obtain
	$$
	\|v\|_{\dot S^1(I_j)} \lesssim  \langle j\eta \rangle  \|f\|_{L_t^1 \dot H^1([0,(j+1)\eta))}
	$$
	Summing in $j$ we obtain \eqref{inhomo}.
\end{proof}

We are now ready to state the following theorem which we will prove by analyzing the linearized equation \eqref{eqnv} around the ground state $W$. 

\begin{theorem}\label{exiuni} 
%New version: 
%Let $c$ be a small constant depending only on the equation \eqref{eqnv}. Let $\lambda\in(0, e_0]$ and $\delta\in(0,\delta_\lambda)$ with $\delta_\lambda:=c\min(\lambda,\lambda^4)$. Then for any $|y_0^-|<\delta$, there exists 
There exists $C>0$ depending only on equation \eqref{eqnv} such that, for any $\lambda\in(0, e_0]$ and $y_0^- \in (-\delta, \delta)$ where 
% $\delta\in(0,\delta_\lambda)$ with 
$\delta=\frac 1C \min(\lambda,\lambda^4)$, there exists 
%
%Let $\lambda\in(0,e_0]$ and $c$ be a sufficiently small constant. 
% Then for any $|y_0^-|<\delta:=c\lambda$, there exists 
 a unique solution to \eqref{eqnv}:
\begin{align*}
v_t=JL v+R(v)
\end{align*}
satisfying 
\begin{equation} \label{E:class}
v(0)=y_0^-V^-+y_0^+V^++v^c(0), \textit{ and } \|v(t) \|_{\dha}
%\dot S^1([t,\infty)\times\R^3)}
\le C \delta e^{-\lambda t}.
\end{equation}
Moreover, 
\begin{align}\label{438}
\begin{cases}
\|v(t)\|_{2}\le C\delta e^{-\lambda t}, \; |y_0^+|+\|v^c(0)\|_{\dha}\le C |y_0^-|^2, \\
\|v\|_{\dot S^1([t,\infty))}\le C^2 \delta e^{-e_0 t}. 
\end{cases}
\end{align}
For any $y_0$, $\tilde y_0$ such that $y_0^-\tilde y_0^->0$ and $|y_0^-|, |\tilde y_0^-|< \delta$, the corresponding solutions $v(t,x)$ and $\tilde v(t,x)$ obey $v(t)=\tilde v(t+T)$ for some $T=T(y_0, \tilde y_0)$. 
\end{theorem}

\begin{proof} 
As from Proposition \ref{prop:linzeng}, $(\dha)^2=E^s\oplus E^u\oplus E^c$, we can decompose 
\begin{align}\label{decu}
v=y^+V^++y^-V^-+v^c
\end{align}
with $y^{\pm}=\langle LV^{\mp}, v\rangle$ and $v^c=v-y^-V^--y^+V^+$. Using the invariance of $JL$ on $E^{u}$, $E^s$ and $E^c$, we reduce the problem to the following system 
\begin{align*}
\begin{cases}
\dot y^-&=-e_0 y^-+R^-(v)\\
\dot y^+&=e_0 y^++R^+(v)\\
\frac{\partial}{\partial t}v^c&=JL v^c+R^c(v). 
\end{cases}
\end{align*}
Here, $R^{\pm}(v)$ and $R^c(v)$ are defined similarly as $y^{\pm}$ and $v^c$. Due to the lack of  exponential decay in the unstable and center directions of the linear flow $e^{tJL}$ as $t \to +\infty$, by Duhamel, exponential decaying solutions must satisfy  
\begin{align}\label{inteqn}
\begin{cases}
y^-(t)&=e^{-e_0 t} y_0^-+\int_0^t e^{-e_0(t-s) }R^-(v(s))ds\\
y^+(t)&=-\int_t^\infty e^{e_0(t-s)}R^+(v(s))ds\\
v^c(t)&=\int_t^\infty e^{JL(t-s)}R^c(v(s))ds. 
\end{cases}
\end{align}

Our goal is to show that the above right sides define a contraction 
\[
(\tilde y^\pm, \tilde v^c) = F(y^\pm, v^c)
\]
on the ball defined by 
\begin{align*}
B_{\delta,\lambda}=\{(y^\pm, v^c) \in C^0([0, \infty)) \times \dot S^1 ([0, \infty)) \mid & \sup_{t\ge 0}e^{\lambda t}|y^{\pm}(t)|\le 2\delta; \\ 
&\sup_{t\ge 0}e^{\lambda t}\|v^c\|_{\dot S^1([t,\infty))}\le 2\delta\}.
\end{align*}
It is easy to see $B_{\delta, \lambda}$ increases in $\delta$ and decreases in $\lambda$. 

 We define another ball $\tilde B=B_{\delta, \lambda}\cap \{v(t,x):\;\sup_{t\ge 0} e^{\lambda t} \|v(t)\|_2\le 2\delta\}$. We will show later that the solution obtained in $B_{\delta, \lambda}$ also belongs to $\tilde B$, from which we immediately prove the $L^2$ regularity in \eqref{438}. 

Taking $(y^{\pm}, v^c)$ from $B_{\delta, \lambda}$, we first reproduce the same bounds on $F(y^\pm, v^c)$ by using the equations \eqref{inteqn}.

To estimate $\tilde y^-(t)$, we first recall that $V^-$ is the eigenfunction of $JL$ associated to the eigenvalue $-e_0<0$, 
%We have 
%$$
%LV^- =e_0 J V^-,
%$$
which allows us to estimate
\begin{align*}
	&|R^-(v(s))| = \left |{\langle LV^+, R(v(s))\rangle } \right | 
	             = \left | {\langle - e_0 JV^+, R(v(s))\rangle}\right | \\
		   \lesssim &\|V^+\|_{6}\|R(v(s))\|_{{\frac 65}} 
		   \lesssim  \|W\|_{6}^3\|v\|_{6}^2+\|v\|_{6}^5\\
		   \lesssim & |y^+(s)|^2+|y^-(s)|^2+\|v^c(s)\|_{\dot H^1}^2 +
			     |y^+(s)|^5+|y^-(s)|^5+\|v^c(s)\|_{\dot H^1}^5 \\
                   \lesssim & (2\delta)^2 e^{-2\lambda s}.
\end{align*}
Inserting this to the first equation in \eqref{inteqn} we have
\begin{align*}
	e^{\lambda t}|\tilde y^-(t)| &\le e^{(\lambda-{e_0})t}|y_0^-| + e^{\lambda t}
	  \int_0^t e^{-(t-s)e_0}|R^-(v(s))|ds \\
	&\le |y_0^-| +C(2\delta)^2\int_0^t e^{(\lambda- e_0)(t-s)}e^{-2\lambda s}ds \\
	&\le |y_0^-| +4C \delta^2/\lambda \le 2 \delta.
\end{align*}
%by using $\delta\le \delta_\lambda$. 
%%if we choose 
%%\begin{align}\label{chd}
%%\delta\le \frac{\lambda}{8C}=:c\lambda. 
%%\end{align}

The estimate of $\tilde y^+(t)$ is similar. Indeed, arguing in the same way as for $R^-(v(s))$, we have 
\begin{align}
	|R^+(v(s))| &\lesssim \delta^2e^{-2\lambda s}, \notag\\
	e^{\lambda t}|\tilde y^+(t)|&\lesssim \delta^2\int_t^{\infty}e^{(\lambda+e_0)(t-s)}e^{-2\lambda s}ds 	\lesssim \delta^2 e^{-\lambda t}
	\le 2\delta\label{1009}
\end{align}
for the same choice of $\delta$. 

We now turn to the estimate of $\tilde v^c(t,x)$ and we start by stating a nonlinear estimate which will be used multiple times.

\begin{claim} \label{rc}For $v$ defined in \eqref{decu} and $\{ y^{\pm}, v^c\} \in B_{\delta,\lambda}$, we have
$$
\|R^c(v(s))\|_{L_t^1\dot H^1([T_0, T_0+T_1])} \lesssim \delta^2 e^{-2\lambda T_0} \langle T_1 \rangle.
$$
\end{claim}

Indeed, from the expression of $v$ in \eqref{decu}, it is straightforward to check 
\begin{align}\label{1208}
\|\nabla v\|_{L_t^5L_x^{\frac {30}{11}}([T_0, T_0+T_1])} \lesssim 
%\langle T \rangle ^{\frac 15}
\delta e^{-\lambda T_0}.
\end{align}
Applying this estimate and using Sobolev embedding, we immediately get
\begin{align}
	\|R(v(s))\|_{L_t^1\dot H^1([T_0, T_0+T_1])}  
	&\lesssim \sum_{i=0}^3 \|W^i v^{5-i}\|_{L_t^1\dot H^1([T_0, T_0+T_1])}\notag\\
		&\lesssim \sum_{i=0}^3 T_1^{\frac i5}\|\nabla W\|_{{\frac {30}{11}}}^i 
	\|\nabla v\|_{L_t^5L_x^{\frac {30}{11}}([T_0, T_0+T_1])}^{5-i} \notag\\
		&\lesssim \sum_{i=0}^3 T_1^{\frac i5}\|\nabla v\|_{L_t^5L_x^{\frac {30}{11}}([T_0, T_0+T_1])}^{5-i}. \label{340}
\end{align}
Inserting \eqref{1208} into \eqref{340}, we proved the Claim \ref{rc}.

We are ready to estimate $\tilde v^c$ on $[T,\infty)$. 
%Without loss of generality, we assume $T\ge 1$. The case when $T<1$ follows from the similar argument with minor changes. 
By triangle inequality, we have
\begin{align*}
	\|\tilde v^c\|_{\dot S^1([T,\infty))}&= \biggl\|\int_t^{\infty} e^{(t-s)JL}R^c(v(s))ds\biggr \|_{\dot S^1([T,\infty))} \\
	%&\le \sum_{N\ge 1, N\in\mathbb N
	%			}\biggl\|\int_t^{\infty} e^{(t-s)JL} R^c(v(s))ds \biggr\|_{\dot S^1([NT, (N+1)T])} \\
	&\le \sum_{
		N\ge 1, N\in\mathbb N
		} \biggl\|\int_t^{T +N}e^{(t-s)JL} R^c(v(s))ds\biggr\|_{\dot S^1([T+ N-1, T+ N])}  \\
	&\quad +\sum_{
		N\ge 1, N\in\mathbb N
		} \biggl\|\int_{T+N}^{\infty} e^{(t-s)JL}R^c(v(s))ds\biggr\|_{\dot S^1([T+N-1, T+ N])} \\
	&:= I + II.
\end{align*}
To estimate $I$, we use time translated version of \eqref{inhomo} in Lemma \ref{lemma3.2} and get 
\begin{align*}
	I &\le \sum_{
		N\ge 1, N\in\mathbb N} 
		%\langle T\rangle^2 
		\|R^c(v(s))\|_{L_t^1\dot H^1([T+N-1, T+N])} \\
	 &\lsm \sum_{
		N\ge 1, N\in\mathbb N		} 
		%\langle T\rangle^3 
		\delta^2 e^{-2\lambda (T+N-1)}  
	 \le \frac 1{\lambda} C \delta^2  e^{-2\lambda T}.
	 % \le \delta e^{-\lambda T}.
\end{align*}
To estimate $II$, we further partition the integral into 
\begin{align*}
	II &\le \sum_{
		N\ge 1	}\sum_{
			 M\ge N+1
		 }\biggl \|\int_{T+M-1}^{T+M} e^{(t-s)JL} 
	     R^c(v(s))ds\biggr\|_{\dot S^1([T+N-1, T+N])} \\
	   &\le \sum_{
			   N\ge 1 ,
			   M\ge N+1
               }\int_{T+M-1}^{T+M} \|e^{(t-s)JL} R^c(v(s))\|_{\dot S^1([T+N-1, T+N])}ds. 
\end{align*}
Note $|t-s|\le M$, applying Lemma \ref{Lemma3.1} we obtain
\begin{align*}               
	   II&\le \sum_{
	N\ge 1, M\ge N+1        
	      }\int_{T+M-1}^{T+M} M^2 \|R^c(v(s))\|_{\dot H^1} ds,
\end{align*}
from which we sum in $N$ and use Claim \ref{rc} to continue	      
	   \begin{align*}   
	   II&\le \sum_{
			   M\ge2
                         }M^3 \|R^c(v(s))\|_{L_t^1\dot H^1([T+M-1, T+M])} \\
	   &\lsm \sum_{
			   M\ge 2
                   	      }M^3 \delta^2 e^{-2\lambda (T+M-1)} 
%	   &\lesssim (2\delta)^2 T^2 e^{-\frac 32 e_0 T} \\
	   \le \frac C{\lambda^4} \delta^2 e^{-2\lambda T}.
	   \end{align*}
Collecting the estimates for $I$ and $II$, 
%and combining the case $T<1$, 
we obtain
\begin{align}\label{1016}
\sup_{T\ge 0}e^{\lambda T}\|\tilde v_c\|_{\dot S^1([T,\infty))}\le \sup_{T\ge 0}e^{\lambda T}(I+II)\le C (\frac 1{\lambda^4} +\frac 1\lambda) \delta^2\le 2\delta. 
\end{align}
%for $\delta\le \delta_\lambda$.

This shows that the map $(\tilde y^\pm, \tilde v^c)$ defined by the right side of \eqref{inteqn} maps $B_{\delta,\lambda}$ to itself. Due to the polynomial form of the nonlinearity, following the similar argument we can easily show the map is a contraction on $B_{\delta,\lambda}$ with a Lipschitz constant $\frac 12$,  hence the existence and uniqueness of the solution to \eqref{inteqn} in $B_{\delta,\lambda}$ is proved.

 Next we show that $(\tilde y^{\pm}, \tilde v^c)\in \tilde B$ if $(y^{\pm},v^c)\in \tilde B$ and it suffices to estimate the $L^2$ norm only. Using the estimate of $\tilde y^{\pm}$ and the fact that $V^{\pm}\in L^2(\R^3) $ from Remark \ref{R:V}, we further reduce the matter to showing $\|\tilde v^c(t)\|_2 \lsm \delta e^{-\lambda t}$. Taking the $L^2$ norm on the expression of $\tilde v^c$ and using the $L^2$ linear estimate from \eqref{607} we have 
\begin{align*}
\|\tilde v^c(t)\|_2\lsm \int_t^\infty \|R^c v(s)\|_2 ds +\int_t^\infty \left\langle t-s\right\rangle^2 \|R^c(v(s))\|_{\dot H^1} ds. 
\end{align*}
From here we partition the integral into pieces and arguing in the same way as above. The only missing piece is $\|R^c v(s)\|_{L_t^1L_x^2}$ on a unit time interval which can be done easily
\begin{align*}
\|R^c(v)\|_{L_t^1L_x^2([t+N-1,t+N])}&\lsm \sum_{i=0}^3\|W\|_{10}^i \|v\|_{L_t^{10} L_x^{10}([t+N-1, t+N])}^{5-i}\\
&\lsm \|v\|_{\dot S^1([t+N-1, t+N])}^2+\|v\|_{\dot S^1([t+N-1, t+N])}^5.
\end{align*}
The rest of the argument will be similar, we omit the details.  This proves the $L^2$ estimate in \eqref{438}.

To see the quadratic estimate \eqref{438}, we note by repeating the same argument, the solution map is contractive on a smaller ball $B_{|y_0^-|,e_0}$. This together with the uniqueness in $B_{\delta,\lambda}$ implies that the constructed solution must lie in $B_{|y_0^-|,e_0}$. From here we apply the estimate in analogue with \eqref{1009} and \eqref{1016} with $\delta$ being replaced by $|y_0^-|$ and $\lambda$ by $e_0$, we immediately obtain 
\begin{align*}
|y_0^+|+\|v^c(0)\|_{\dha}\lsm |y_0^-|^2. %\lsm \frac 1{e_0^4} |y_0^-|^2
\end{align*}

In the above we prove the existence of the stable solution with $y^-(0)=y_0^-$ which is unique in $B_{\delta, \lambda}$. The stronger uniqueness of such solution in the set of functions characterized by \eqref{E:class} is a simple consequence of the following Lemma \ref{L:decay} and the above uniqueness. 

To complete the proof of Theorem \ref{exiuni}, for any $y_0^-$ and $\tilde y_0^-$ satisfying $|y_0^-|, |\tilde y_0^-|< \delta$ and $y_0^-\tilde y_0^->0$, let $v(t)$ and $\tilde v(t)$ be the corresponding exponentially decaying solutions. From the continuity and decay of $\tilde v(t)$ in $\dha(\mathbb{R}^3)$, we know there must exist a time $T$ such that $\tilde y^-(T)=\langle LV^+, \tilde v(T)\rangle= y_0^-$, from the uniqueness we conclude that $v(t)=\tilde v(t+T)$. 
\end{proof}

The following lemma gives the exponential decay of the Strichartz norm from the exponential decay of the $\dha$ norm.

\begin{lemma} \label{L:decay}
Assume $v(t)$ is a solution to \eqref{eqnv} satisfying that, for some $\lambda>0$,   
\[
\|v(t) \|_{\dha} \lsm e^{-\lambda t}, \quad t\ge 0,
\]
then
\begin{align}\label{101}
\| v\|_{\dot S^1 ([t, \infty))} \lsm e^{-\lambda t}, \quad t\ge 0.
\end{align}
\end{lemma}

\begin{proof}
Let $T_0>0$ be sufficiently large. It suffices to prove the estimate \eqref{101} for all $t\ge T_0$ as the estimate for $t\in [0,T_0)$ follows from the estimate of $\|v\|_{\dot S^1([T_0,\infty))}$ and the standard local estimate on $[t,T_0)$ . Let $\tau\ge T_0$ and $\eta$ be a small number to be chosen later. Applying the Strichartz estimate on the interval $[\tau, \tau+\eta]$ and using the similar nonlinear estimate as in \eqref{340}, we obtain 
\begin{align*}
\|v\|_{\dot S^1{([\tau,\tau+\eta])}}\le C\|v(\tau)\|_{\dha}+C\sum_{i=0}^4 \eta^{\frac i5} \|v\|_{\dot S^1{([\tau,\tau+\eta])}}^{5-i}
\end{align*}
for some constant $C$ independent of $v$ and $\eta$. Recall that $\|v(\tau)\|_{\dha}\lsm e^{-\lambda \tau}\le e^{-\lambda T_0}$, for $\eta$ sufficiently small and $T_0$ sufficiently large, the standard continuity argument gives 
\begin{align*}
\|v\|_{\dot S^1([\tau,\tau+\eta])}\le 2C\|v(\tau)\|_{\dha}\lsm e^{-\lambda \tau}. 
\end{align*}
The estimate of $\|v\|_{\dot S^1([t,\infty))}$ then comes from partitioning the interval $[t,\infty)$ and adding up the estimate on each subinterval. The proof is complete. 

%
%
%Without loss of generality, we may assume that $\| v\|_{L_t^\infty \dha} \le \delta$ for some small $\delta$ to be determined in the below. For any $t_0 \in (0, 1]$ and $t\ge t_0$, according to Lemma \ref{Lemma3.1} and \ref{lemma3.2}, we have 
%\begin{align*}
%\|v\|_{\dot S^1 ([t-t_0, t+t_0])} \lsm& \| v(t)\|_{\dot H^1} + \|R(v)\|_{L_t^1 \dot H^1 ([t-t_0, t+t_0])} \\
%\lsm & \| v(t)\|_{\dot H^1}  + \|v\|_{\dot S^1 ([t-t_0, t+t_0])}^2 + \|v\|_{\dot S^1 ([t-t_0, t+t_0])}^5.  
%\end{align*}
%where the constant is independent of $\delta$ and $t$ and we used the type of argument proving Claim \ref{rc} to estimate the $R(v)$ term. When the upper bound $\delta>0$ of $\| v\|_{L_t^\infty \dha}$ is sufficiently small, the above estimate implies 
%\[
%\|v\|_{\dot S^1 ([t-1, t+1])} \lsm \| v(t)\|_{\dot H^1} \lsm e^{-\lambda t}, \quad \forall t\ge 1. 
%\] 
%Summing up this estimate yields
%\[
%\| v\|_{\dot S^1 ([t, \infty))} \le \sum_{n=0}^\infty \| v\|_{\dot S^1 ([t+ 2n, t+2(n+1)])} \lsm e^{-\lambda t}
%\]
%and the proof is complete. 
\end{proof}

Lemma \ref{L:decay} together with Theorem \ref{exiuni} finally gives rise to the following result, which characterizes all solutions decaying exponentially to the ground state:

\begin{corollary}\label{cor:411}
There exist exactly two solutions (up to time translation) $W_{\pm}$ of $NLS_a$ satisfying 
\begin{align*}
\begin{cases}
\|W_{\pm}-W\|_{H^1}\le C e^{-e_0t}, \ \forall t\ge 0. \\
\|W_+(0)\|_{\dha}>\|W\|_{\dha},\ \ \|W_-(0)\|_{\dha}<\|W\|_{\dha}.
\end{cases}
\end{align*}
Moreover, if a solution $u(t,x)$ of $NLS_a$ satisfying 
\begin{align*}
\|u(t)-W\|_{\dha}\le C e^{-\lambda t},\ \forall t\ge 0
\end{align*}
for any $C>0$ and $\lambda\in (0, e_0]$, there must exist unique $T^{\pm}$ such that
\begin{align*}
\begin{cases}
u(t)=W_+(t+T^+) & \textit{ if } \|u\|_{\dha}>\|W\|_{\dha},\\
u(t)=W_-(t+T^-) & \textit{ if } \|u\|_{\dha}<\|W\|_{\dha}. 
\end{cases}
\end{align*}
\end{corollary}

We remark that this Corollary does not tell us the behavior of $W_{\pm}$ for $t<0$, we will discuss this problem in Section 7 and Section 8, and complete the picture of the dynamics of all solutions on the energy surface.

\section{Global analysis-Virial}

In the previous sections, we develop the modulation analysis which enables us to control the solution near the two dimensional manifold generated by the symmetry transformations applied to $W$. When the solution is away from the manifold, we use the monotonicity formula arising from Virial to control the solution. To this end, in this section we establish Virial estimates by incorporating the modulation estimates developed in Section 4.

Let $\phi(x)$ be a smooth radial function such that 
\begin{equation*}
\phi (x)=%
\begin{cases}
|x|^{2}, & |x|\leq 1; \\ 
0, & |x|>2,%
\end{cases}%
\;  \textit{ and }\phi _{R}(x)=R^{2}\phi \bigl(\frac{x}{R}\bigr).
\end{equation*}%
 Moreover, we can choose $\phi$ such that the radial derivative satisfies
\begin{align}\label{510}
\phi''(r)\le 2. 
\end{align}
From such $\phi$ we define the truncated Virial 
\begin{equation*}
V_{R}(t)=\int_{\mathbb{R}^{3}}\phi _{R}(x)|u(t,x)|^{2}dx.
\end{equation*}%
For a solution $u(t)$ of NLS$_a$ with $E(u)=E(W)$, the time derivatives of $V_R(t)$ are
computed as  
\begin{align}
\partial _{t}V_{R}(t)& =2\mathbf{Im}\int_{\mathbb{R}^{3}}\overline{u(t)}%
\,\nabla u(t)\cdot \nabla \phi _{R}dx;  \notag \\
\partial _{tt}V_{R}(t)& =4\mathbf{Re}\int_{\mathbb{R}^{3}}(\phi
_{R})_{jk}(x)u_{j}(t)\bar{u}_{k}(t)\,dx-\frac{4}{3}\int_{\mathbb{R}%
^{3}}(\Delta \phi _{R})|u(t)|^{6}\,dx  \notag \\
& \quad -\int_{\mathbb{R}^{3}}(\Delta^2 \phi
_{R})|u(t)|^{2}\,dx+4a\int_{\mathbb{R}^{3}}\tfrac{x}{|x|^{4}}\nabla \phi
_{R}|u(t)|^{2}\,dx  \notag \\
& =16(\Vert W\Vert _{\dot{H}_{a}^{1}}^{2}-\Vert u(t)\Vert _{\dot{%
H}_{a}^{1}}^{2})+A_{R}(u(t)), \notag
\end{align}%
where 
\begin{eqnarray*}
A_{R}(u(t)) &=&\int_{|x|>R}(4\partial_{rr} \phi_R-8)|\nabla
u(t)|^{2}dx+\int_{|x|>R}(-\frac{4}{3}\Delta \phi _{R}+8)|u(t)|^{6}dx \\
&&-\int_{\mathbb{R}^{3}}\Delta^2\phi _{R}|u(t)|^{2}dx+\int_{|x|>R}(%
\frac{4a}{|x|^{4}}x\nabla \phi _{R}|u(t)|^{2}-\frac{8a|u(t)|^{2}}{|x|^{2}}%
)dx.
\end{eqnarray*}%
As seen in Section 4, $\bd(u(t))$ plays a role of measuring the distance between $u(t)$ and the manifold, we then rewrite $\partial_{tt}V_R(t)$ into  
\begin{align}\label{1122}
\partial_{tt}V_R(t)=\begin{cases}
16 \ \bd(u(t))+A_R(u(t)),& \textit{ if } \| u(t)\|_{\dha}<\| W\|_{\dha};\\
-16\  \bd(u(t))+A_R(u(t)), & \textit{ if } \| u(t)\|_{\dha}>\| W\|_{\dha}.
\end{cases}
\end{align}

The rest of this section is devoted to giving proper estimates on $\partial_t V_R(t)$ and $A_R(t)$. We start with the following elementary lemma which shows how they are rescaled under the transformation of symmetries. 

\begin{lemma}
\label{lem:1103} For any $(\theta,\mu)\in \mathbb{S}^1\times\R^+$, we have the following scaling relations: 
%alternative expressions 
\begin{align*}
\partial _{t}V_{R}(t)& =\mu ^{-2}\cdot 2\mathbf{Im}\int_{\mathbb{R}%
^{3}}\nabla \phi _{\mu R}\cdot \nabla u_{[\theta ,\mu ]}\ \overline{%
u_{[\theta ,\mu ]}}dx. \\
A_{R}(u(t))& =A_{\mu R}(u(t)_{[\theta ,\mu ]}).
\end{align*}
In addition, 
$$
A_{R}(W)=0.
$$
\end{lemma}

The verification of this Lemma is straightforward so we skip it. This Lemma together with the modulation analysis from Section 4 yields:

\begin{lemma}[Virial estimate]
\label{lem:1104} Let $u(t)$ be an $\dha$-solution of $NLS_a$ with $E(u)=E(W)$. For those $t$ satisfying $\bd(u(t))<\delta_0$, let 
$$
u(t)_{[\smt]}=W+v(t)
$$
be the orthogonal decomposition of $u(t)$ given by Lemma \ref{L: Modulation} with the corresponding bounds. We have 

\begin{align}
|\partial _{t}V_{R}(t)|&\lesssim R^{2}\bd(u(t)) ,\label{gap 1}
\end{align}
 \begin{align}\label{gap 3}
A_R(u(t))\lsm \int_{|x|>R}\bigr(|u(t,x)|^6
+\frac{|u(t,x)|^{2}}{|x|^{2}}\bigr)dx
\end{align}
\begin{align}
|A_{R}(u(t))|&\lesssim\begin{cases} \int_{|x|>R}\bigl(|\nabla u(t)|^{2}+|u(t)|^6
+\frac{|u(t)|^{2}}{|x|^{2}}\bigr)dx.\\
[(\mu (t)R)^{-\frac{\beta}{2}}\bd(u(t))+%
\bd(u(t))^{2}],\ \ \textit{ if  }\bd(u(t))<\delta _{0} \textit{ and } |\mu(t)R|\gtrsim 1, 
\end{cases}\label{gap 2}
\end{align}
where the constants  are independent of $\bd(u)$, $R$ and $\|u\|_{\dha}$. 
\end{lemma}

\begin{proof}
We first estimate $\partial _{t}V_{R}(t)$. From H\"{o}lder inequality and Sobolev embedding, we have
\begin{equation*}
|\partial _{t}V_{R}(t)|\leq \Vert \nabla u\Vert _{2}\Vert u\Vert _{6}\Vert \nabla \phi _{R}\Vert _{3}\lesssim R^{2}\Vert \nabla \phi \Vert
_{3}\Vert \nabla u\Vert _{2}^{2}\lesssim R^{2}(\dut+\|W\|_{\dha}^2).
\end{equation*}%
This proves the bound in the case of $\dut\ge \delta_0$ if the implicit constant is allowed to depend on $\delta_0$. To get the bound in the rest of the case $\dut<\delta_0$, we use Lemma \ref{lem:1103} with $(\theta,\mu)$ being given by $(\theta(t), \mu(t))$ and the decomposition to get 
\begin{align*}
|\partial _{t}V_{R}(t)|& =\bigl|\mu (t)^{-2}\cdot 2\mathbf{Im}\int_{\mathbb{R}^{3}}\nabla \phi _{\mu (t)R}\nabla (W+v(t))(W+\bar v(t))dx\bigr| \\
& =\bigl|\mu (t)^{-2}\cdot 2\mathbf{Im}\int_{\mathbb{R}^{3}}\nabla \phi
_{\mu (t)R}(\nabla W\bar{v}(t)+\nabla v(t)W+\nabla v(t)\bar{v}(t))dx\bigr| \\
& \leq \mu (t)^{-2}\Vert \nabla \phi _{\mu (t)R}\Vert _{3}(\Vert \nabla
W\Vert _{2}\Vert \nabla v\Vert _{2}+\Vert \nabla v\Vert _{2}^{2}) \\
& \lesssim R^{2}\mathbf{d}(u(t)).
\end{align*}%

We now turn to estimating $A_{R}(u(t))$. Using \eqref{510}, we can throw away the first non-positive term in the expression of $A_R(u(t))$ and estimate the rest three terms to get \eqref{gap 3}. The same direct estimate also gives the first line in \eqref{gap 2}. To get the second bound when $\mathbf{d}(u(t))<\delta
_{0}$ and $\mu (t)R\gtrsim 1$, we recall $W(x)=O(|x|^{-\frac 12-\frac 12\beta})$ for $|x|\gtrsim 1$. This together with the decomposition and Lemma \ref{lem:1103} yields
\begin{align*}
|A_{R}(u(t))|& =|A_{\mu (t)R}(u(t)_{[\theta (t),\mu (t)]})|=|A_{\mu
(t)R}(W+v(t))-A_{\mu (t)R}(W)| \\
& \lesssim \Vert \nabla W\Vert _{L^{2}(|x|\geq \mu (t)R)}\Vert \nabla
v(t)\Vert _{2}+\Vert \nabla v(t)\Vert _{2}^{2} \\
& +\Vert v\Vert _{6}(\Vert W\Vert _{L^{6}(|x|\geq \mu
(t)R)}^{5}+\Vert v(t)\Vert _{6}^{5}%
)+\Vert {W}/{|x|}\Vert _{L^{2}(|x|\geq \mu (t)R)}\Vert \nabla v(t)\Vert
_{2} \\
& \lesssim (\mu (t)R)^{-\frac \beta {2}}\mathbf{d}(u(t))+(\mathbf{d}%
(u(t)))^{2}.
\end{align*}%
Lemma \ref{lem:1104} is proved.
\end{proof}

\section{Exponential convergence in the sub-critical case\label{decay_sub}}

In this section, we focus on characterizing the non-scattering solutions on the energy surface of $E(W)$ when the kinetic energy is less than that of the ground state $W$. 
The main result is the following

\begin{theorem}
\label{thm: subcase} Let $u$ be a solution of NLS$_{a}$ satisfying %
\begin{equation}
E(u)=E(W),\Vert u_{0}\Vert _{\dot{H}_{a}^{1}}<\Vert W\Vert _{\dot{H}%
_{a}^{1}},\Vert u\Vert _{S([0,\infty ))}=\infty.  \label{cond}
\end{equation}%
Then there exist $\theta \in \mathbb S^1$, $\mu >0$ and a unique time $T=T(u)$ such that 
\begin{equation}\label{357}
u(t,x)= e^{i\theta}\mu^{\frac 12}W_-(\mu^2 t+T,\mu x).
\end{equation}
In the opposite time direction, $u$ exists globally and obeys $\|u\|_{S((-\infty,0])}<\infty$.
%Moreover, there exists a unique time $T=T(u)$ such that 
\end{theorem}

We note first that \eqref{357} comes directly from 
\begin{equation}\label{448}
\Vert u(\mu^{-2}t)_{[\theta ,\mu ]}-W\Vert _{\dot{H}_{a}^{1}}\leq Ce^{-ct},\forall \
t\geq 0,
\end{equation}
%\eqref{448} 
satisfied by the solution $u(\mu^{-2} t)_{[\theta, \mu]}$ and Corollary \ref{cor:411}. 
%Indeed, assuming \eqref{448}, we know that the solution 
%\[\mathbf{g}_{\theta, \mu}u(t,x)=e^{i\theta}\mu^{-\frac 12}u(t/\mu^2,x/\mu)\] satisfies 
%\[\|\mathbf{g}_{\theta,\mu}u(t)-W\|_{\dha}\le C e^{-c/\mu^2 t}.\]
%Applying \eqref{cor:411}, we know immediately that there exist a unique time $T$ such that
%\[\mathbf{g}_{\theta, \mu}u(t,x)=W_-(t+T,x)\] 
%which after undoing the scaling, verifies the last equation in Theorem \ref{thm: subcase}. 

Therefore throughout the rest of this section we will only focus on the proof of \eqref{448}.
We start by discussing properties 
of solutions obeying \eqref{cond}.

\subsection{Properties of solutions satisfying \eqref{cond}}

From \eqref{cond}, we know $u$ is non-scattering at the minimal energy $E(W)$.  The minimality induces the compactness at least in the radial case, as was proved in an earlier work \cite{R: Monica inverse energy, R:yk_general, R:yk_radial}. In the non-radial case and dimension $d=3$, the compactness is still unavailable we will take \eqref{precomp} as an assumption and build our conditional result upon it. Results in dimensions four and five become unconditional. More specifically, 
there exists $\lambda (t):[0,\infty )\rightarrow \mathbb{R}^{+}$ such that 
\begin{equation}
\{u(t)_{[\lambda (t)]}\}_{t\in \lbrack 0,\infty )}\textit{ is precompact in }%
\dha(\R^3)(\textit{or } \dot H^1(\R^3)).  \label{precomp}
\end{equation}%

The first step of proving such statement is to take an arbitrary sequence $\{u(t_n)\}$ and show that there exist $\lambda_n$ such that $\{\bg_{\lambda_n}u(t_n)\}$ is precompact in $\dot H^1(\mathbb{R}^3)$. While this had been achieved in \cite{R: Monica inverse energy,R:yk_general, R:yk_radial}, it is not entirely clear from here how to jump to the continuous choice of $\lambda(t)$. Here we provide a point of view through which we are able to make the choice of continuous $\lambda(t)$ explicitly and more quantitatively. 

Let $\psi: [0, \infty) \to [0, 1)$ be a smooth function  such that 
\[
\psi(0)=0, \; \lim_{s\to +\infty} \psi(s)= 1 \; \textit{ and }\; \psi'(s)>0. 
\]
Define the weighted norm \[V(R, u)=\int_{\R^3}\psi(|x|/R)|\nabla u|^2 dx.\] Then for any nontrivial function $u\in \dot H^1(\mathbb{R}^3)$, we can easily check that
\begin{align*}
\partial_R V(R, u)<0,\ V(0,u)=\int_{\R^3}|\nabla u|^2 dx,\\
 V(\infty, u)=0, \mbox{ and } V(R, u_{[\mu]})=V(R/\mu, u).
 \end{align*}
Due to the monotonicity of $V$, for any $u\in \dot H^1(\mathbb{R}^3)$, there exists a unique $\Lambda(u)$ such that $V(1, u_{[\Lambda]})=\frac 12$. Clearly $\Lambda: \dot H^1(\mathbb{R}^3) \to \R^+$ is smooth. We have the following lemma. 

\begin{lemma} 
Suppose sequences $\{u_n\} \subset \dot H^1(\R^3)$ and $\{\lambda_n\} \subset \R^+$ satisfy that $\bg_{\lambda_n}u_{n}$ converges in $\dot H^1(\R^3)$. Then $\{\bg_{\Lambda(u_{n})}u_{n}\}$ also converges in $\dot H^1(\R^3)$. 
\end{lemma}

\begin{proof} 
Let
\begin{align}\label{dis1}
\lim_{n\to \infty} \bg_{\lambda_n}u_n=\phi \textit{ in } \dot H^1(\R^3).
\end{align}
Let $\Lambda_0=\Lambda(\phi)$. Then from scaling we have
\begin{align*}
V(1/\Lambda_0,\phi)=V(1, \bg_{\Lambda_0}\phi)=V(1,\bg_{\Lambda(u_n)} u_n)=V({\lambda_n}/{\Lambda(u_n)},\bg_{\lambda_n }u_n)
\end{align*}
which clearly implies $\lambda_n/\Lambda(u_n)\to 1/\Lambda_0$ as $n\to \infty$ by using the strong convergence \eqref{dis1}. Hence 
\[\bg_{\Lambda(u_n)}u_n\to \bg_{\Lambda_0}\phi \mbox{ in } \dot H^1(\R^3).\]
%The precompactness of $\{\bg_{\mu(t)}u(t)\}$ is proved.
\end{proof}

Therefore for any solution $u(t)$ whose orbit is precompact modular scaling in $\dot H^1(\R^3)$, we can take $\lambda(t)=\Lambda(u(t))$ as an underlying choice of scaling parameter which is subject to further mollification in this section. On the one hand, the precompactness of $u(t)$ up to the rescaling has some crucial implications on this $\lambda(t)$. On the other hand,  there is certain freedom in the choice of this scaling function $\lambda(t)$ and we will refine our choice to help us to prove Theorem \ref{thm: subcase}. 

Firstly the compactness implies directly that there exists $C(\varepsilon )>0$ such that 
\begin{equation}
\int_{|x|>\frac{C(\varepsilon )}{\lambda (t)}}|\nabla u(t)|^{2}+|u(t)|^{6}+\frac{|u(t)|^{2}}{|x|^{2}}dx<\varepsilon .  \label{AA char}
\end{equation}%

%some of the properties of $\lambda(t)$ indicated by the compactness. 
Secondly we recall that the scaling parameter $\lambda (t)$ obeys 
%(see for example \cite{})
\begin{equation}
\lim_{t\rightarrow \infty }\lambda (t)\sqrt{t}=\infty \label{lambda t}.
\end{equation}%
Indeed, if this is not true, there exists a sequence of time $t_n\to \infty$ such that $\lambda^2(t_n)t_n\to c<\infty$, as a result
\begin{equation}\label{1216}
\lim_{n\to \infty}\lambda(t_n)\to 0.
\end{equation}
Let 
\begin{align*}
v_n(t)=\bg_{\lambda(t_n)}u\bigl(t_n+\tfrac t{\lambda^2(t_n)}\bigr),
\end{align*}
we have 
\begin{equation}\label{1217}
v_n(0)=\bg_{\lambda(t_n)}u(t_n), \ v_n(-t_n\lambda^2(t_n))=\bg_{\lambda(t_n)}u(0).
\end{equation}
From the compactness, there exists a subsequence and $v_0\in \dha(\mathbb{R}^3)$ such that $v_n(0)\to v_0$ in $\dha(\mathbb{R}^3)$. Let $v(t,x)$ be the solution of $NLS_a$ with data $v_0$. The standard local theory implies $v_n(-t_n\lambda^2(t_n))\to v(-c)\neq 0$ in $\dha(\mathbb{R}^3)$, which immediately contradicts with $ \bg_{\lambda(t_n)}u(0)\rightharpoonup 0$ weakly in $\dha(\mathbb{R}^3)$ from \eqref{1216}. 
%As seen from \eqref{AA char}, the choice of scaling function $\lambda(t)$ is obviously not unique as it can be mollified within the constant range of $\lambda(t)$. This gives us a great deal of room in choosing the specific candidate to better serve our goal. Our ultimate choice of $\lambda(t)$ comes from the following two observations on $\lambda(t)$. The first one is the almost constancy of $\lambda(t)$. 
Next, we have the following

\begin{lemma}[Almost constancy] \label{lem:almost constancy} 
Let $u$ be the solution satisfying \eqref{cond}, then there exist $\delta>0$ and $0<c<C<\infty$ such that for any $\tau\ge 0$, on the interval
\begin{align*}
I_\tau:=[\tau,\tau+\frac \delta {\lambda^2(\tau)}],
\end{align*}
we have
\begin{align}\label{221}
c\le \frac{\lambda(\tau_1)}{\lambda(\tau_2)}\le C, \textit{ for any } \tau_1,\tau_2\in I_\tau.
\end{align}

%ii) Fix a constant $\delta>0$. Then there exists $c(\delta)>0$ such that if 
%\begin{align*}
%\sup_{t\in I_\tau}\dut >\delta, 
%\end{align*}
%we must have 
%\begin{align*}
%\min_{t\in I_\tau} \dut>c(\delta). 
%\end{align*}
\end{lemma}

\begin{proof} 
%Almost constancy of the scaling parameters has been observed in \cite{}, we give the proof for the sake of completeness. 

We argue by contradiction. Suppose \eqref{221} fails, there must exist two sequences of times $0<t_n<s_n<\infty$ and $(s_n-t_n)\lambda^2(t_n)\to 0$ but
\begin{equation}\label{231}
\frac{\lambda(t_n)}{\lambda(s_n)}+\frac{\lambda(s_n)}{\lambda(t_n)}\to \infty.
\end{equation}
Define the scaled solution
\begin{align}\label{318}
v_n(t,x)=\bg_{\lambda(t_n)}u\bigl(t_n+\frac t{\lambda^2(t_n)}\bigr) \textit{ and } \gamma_n:=\lambda^2(t_n)(s_n-t_n). 
\end{align}
We have 
\begin{align}\label{247}
\begin{cases}
v_n(0)=\bg_{\lambda(t_n)}u(t_n),\\
v_n(\gamma_n)=\bg_{\lambda(t_n)}u(s_n)=\bg_{\lambda(t_n)/\lambda(s_n)}\bg_{\lambda(s_n)}u(s_n),\\
\gamma_n\to 0.
\end{cases}
\end{align}
From the first equation in \eqref{247} and the compactness, we know there exist a subsequence and $v_0\in\dot H_a^1(\mathbb{R}^3)$ such that $v_n(0)\to v_0$ in $\dot H_a^1(\mathbb{R}^3)$. This together with the standard local theory implies
\begin{align*}
\lim_{n\to \infty} v_n(\gamma_n)=v_0\neq 0 \mbox{ in } \dot H_a^1(\mathbb{R}^3).
\end{align*}
In addition, the second expression in \eqref{247} together with  \eqref{231} and the compactness along the sequence $\{s_n\}$ imply 
\begin{align*}
v_n(\gamma_n)\rightharpoonup 0\textit{ weakly in }  \dot H_a^1(\mathbb{R}^3),
\end{align*}
after passing to a subsequence if necessary. We get a contradiction. Lemma \ref{lem:almost constancy} is then proved. 
%The proof of ii) follows along the same line which we sketch here. Suppose statement ii) fails, there exist a sequence of intervals $I_{\tau_n}$ and $t_n,s_n\in I_{\tau_n}$ such that
%\begin{align}\label{316}
%\mathbf{d}(u(s_n))>\delta, \mbox{ but } \lim_{n\to\infty}\mathbf{d}(u(t_n))=0. 
%\end{align}
%Define $v_n$, $\gamma_n$ as in \eqref{318}, we have the same set of equations as in \eqref{247}. Again, from the first equation in \eqref{247}, we have $v_n(0)\to v_0$ in $\dot H_a^1$ thus $\mathbf{d}(v_0)=0$. From the variational characterization of $W$, we know $v_0=W_{[\theta_0,\lambda_0]}$. 
%
%On the other hand, the second equation in \eqref{247} together with the local theory imply
%\begin{align*}
%u(s_n)_{[\lambda(t_n)]}=v_n(s_n)\to W_{[\theta_0,\lambda_0]} \mbox{ in }\dot H_a^1. 
%\end{align*}
%Hence $\mathbf{d}(u(s_n))\to 0$ which contradicts with \eqref{316}. Lemma \ref{lem:almost constancy} is proved. 
\end{proof}

The next observation on $\lambda(t)$ is that $\lambda (t)$ is basically comparable to $\mu(t)$ given by Proposition \ref{P:vcW} when the solution $u(t)$ is close to the manifold. 

\begin{lemma}
\label{cop} 
Let $u$ be the solution of NLS$_a$ on the time interval $I$ satisfying \eqref{precomp}. Suppose $\bd(u(t))<\delta _0$ on $I$ hence $u(t)$ is subject to the orthogonal decomposition $\bg_{\theta(t),\mu(t)}u(t)=W+v(t)$. Then there exist constants $0<c<C<\infty$ such that
\begin{equation*}
c<\frac{\lambda(t)}{\mu(t)}<C, \ \forall t\in I.
\end{equation*}%
%In particular, we can redefine $\lambda(t)$ so that $\lambda(t)=\mu(t)$ on $I$. 
\end{lemma}

\begin{proof}
We argue by contradiction. Suppose this is not true, there must
exist a sequence of times $t_{n}\in I$ such that 
\begin{equation}\label{836}
\frac{\mu (t_{n})}{\lambda (t_{n})}\rightarrow 0\ \textit{ or } \ \frac{\mu (t_{n})%
}{\lambda (t_{n})}\rightarrow \infty .
\end{equation}%
From the compactness we can extract a subsequence and $V\in \dot H_a^1(\mathbb{R}^3)$ such that
\begin{equation*}
\bg_{\lambda(t_n)}u(t_{n})\rightarrow V\mbox{ in }\dot{H}_{a}^{1}(\mathbb{R}^3),
\end{equation*}%
which along with $\bd( u(t))<\delta_0$ implies
\begin{equation}
\Vert V\Vert _{\dot{H}_{a}^{1}}^{2}\longleftarrow \Vert u(t_{n})\Vert _{\dot{%
H}_{a}^{1}}^{2}=\Vert W\Vert _{\dot{H}_{a}^{1}}^{2}-\mathbf{d}%
(u(t_{n}))>\Vert W\Vert _{\dot{H}_{a}^{1}}^{2}-\delta _{0}.  \label{V up}
\end{equation}%
On the other hand, along the same sequence, we apply the symmetry $\bg_{-\theta(t_n),\frac{\lambda(t_n)}{\mu(t_n)}} $ on both sides of the orthogonal decomposition  
\begin{equation}\label{842}
\bg_{\theta (t_{n}),\mu (t_{n})}u(t_{n})=W+v(t_{n}),
\end{equation}
to obtain
\begin{equation*}
\bg_{\lambda(t_n)}u(t_{n})=\bg_{-\theta(t_n),\frac{\lambda(t_n)}{\mu(t_n)}} (W+v(t_{n})).
\end{equation*}%
Passing to a subsequence if necessary and taking weak limit on both sides, using \eqref{836} we have 
\begin{equation*}
\bg_{-\theta(t_n),\frac{\lambda(t_n)}{\mu(t_n)}} v(t_{n})\rightharpoonup V \mbox{ weakly in } \dot{H}_{a}^{1}(\mathbb{R}^3).
\end{equation*}%
 This together with \eqref{modulation 1} shows
\begin{equation*}
\Vert V\Vert _{\dot{H}_{a}^{1}}\leq \liminf_{n\rightarrow \infty }\Vert
\bg_{-\theta(t_n),\frac{\lambda(t_n)}{\mu(t_n)}} v(t_{n})\Vert _{\dot{%
H}_{a}^{1}}\lesssim\bd(u(t_{n}))<\delta _{0},
\end{equation*}%
which contradicts (\ref{V up}) and completes the proof of this lemma.
\end{proof}

Next we show that such precompactness implies that $u(t)$ keeps getting closer to the manifold $\{\bg_{\sm}W\}$.
\begin{lemma}
\label{L: sub tn}
Let $u$ be the solution of NLS$_{a}$ satisfying %
\eqref{cond}. Then there exists a sequence of time $t_{n}\rightarrow \infty $ such
that $\mathbf{d}(u(t_{n}))\rightarrow 0$.
\end{lemma}

\begin{proof}
Let $C(\eps)$ be the function defined in \eqref{AA char}. Then from \eqref{lambda t}, for any $\eps >0$, there exists $T_0=T_{0}(\eps)>0$ such that when $t>T_0$, $\lambda(t)t^{\frac 12}>{C(\eps)}/{\eps^{\frac 12}}$. Therefore on the time interval $[T_0, T]$ we have 
\begin{align*}
(\eps T)^{\frac 12}>\frac{C(\eps)}{\lambda(t)}, \ \forall t\in [T_0, T]. 
\end{align*}
Take $R=\left( \eps T\right) ^{\frac{1}{2}}$ and apply Lemma \ref%
{lem:1104} for $t\in[ T_{0}, T]$ we obtain
\begin{align*}
|\partial _{t}V_{R}(t)|& \lesssim R^{2}=\eps T, \\
|A_{R}(u(t))|& \lesssim \int_{|x|>R}\bigl(|\nabla u(t)|^{2}+|u(t)|^{2^{\ast
}}+\frac{|u(t)|^{2}}{|x|^{2}}\bigr)dx \\
& \lesssim \int_{|x|>\frac{C(\eps )}{\lambda (t)}}\bigl(|\nabla
u(t)|^{2}+|u(t)|^{2^{\ast }}+\frac{|u(t)|^{2}}{|x|^{2}}\bigr)dx \\
& \leq \eps,
\end{align*}%
These two estimates together with \eqref{1122} and \eqref{gap 2} give 
\begin{equation*}
\partial _{tt}V_{R}(t)\geq 16\mathbf{d}(u(t))-C\eps .
\end{equation*}%
Integrating in $t$ over $[T_0,T]$ and dividing by $T$ we have 
\begin{equation*}
\frac{1}{T}\int_{T_{0}}^{T}\mathbf{d}(u(t))dt\lesssim \frac{\eps
(T-T_{0})+R^{2}}{T}\lesssim \eps,
\end{equation*}%
which immediately gives
\begin{equation*}
\lim_{T\rightarrow \infty }\frac{1}{T}\int_{0}^{T}\mathbf{d}(u(t))dt=0,
\end{equation*}%
by first taking $T\to \infty$ then $\eps\to 0$. The convergence of $\bd(u(t))$ along a sequence of time is proved. 
\end{proof}

Lemma \ref{lem:almost constancy} implies that $\lambda(t)$ has a change of $C\lambda(t)$ on the interval with the length $O(1/\lambda^2(t))$. Therefore it is intuitive to imagine 
%we can choose $\lambda(t)$ such that
\begin{align}\label{chosel1}
|\lambda'(t)|\le C\lambda^3(t). 
\end{align}
%We now return to \eqref{precomp} to obtain more properties of $\lambda(t)$. 
Lemma \ref{cop} implies that we can replace $\lambda(t)$ by $\mu(t)$ on the interval where $\bd(u(t))<\delta_0$. From \eqref{chosel1} and the derivative estimate for $\mu(t)$ in Lemma \ref{L: Modulation}, it is reasonable to expect 
%we conclude that $\lambda(t)$ can be taken to satisfy the following
\begin{align}\label{ourl}
\frac{|\lambda'(t)|}{\lambda^3(t)}\le 
\begin{cases}
C, &\mbox{ when } \dut>\delta_0;\\
C\dut,&\mbox{ when } \dut\le \delta_0,
\end{cases}
\end{align}
from which we may further modifying the constant to guarantee 
\begin{align}\label{perpl}
\frac{|\lambda'(t)|}{\lambda^3(t)}\le C\dut, \ \forall t\ge 0.
\end{align}
In fact we can modify $\lambda(t)$ such that it is differentiable almost everywhere and 
\begin{align}\label{intel}
\biggl|\frac 1{\lambda^2(a)}-\frac 1{ \lambda^2(b)}\biggr|\le C\int_a^b \mathbf{d}(u(t))dt, \ \forall [a,b]\subset[0,\infty).
\end{align}
See Lemma \ref{L:lambda-1} in the Appendix. 

We will revisit this estimate later when we prove the uniform lower bound for $\lambda(t)$. Now we turn to considering the distance function $\dut$ with the goal of proving the exponential decay of $\dut$. We start by showing 

\begin{lemma}[Integral estimate of $\dut$]
\label{L: Virial}
Let $u$ be the solution of NLS$_a$ satisfying \eqref{precomp}, then there exists $C>0$ such that for any $[a,b]\subset [0,\infty)$,
\begin{equation}\label{1019}
\int_{a}^{b}\mathbf{d}(u(t))dt\leq C\sup_{t\in
[a,b]}\frac{1}{\lambda (t)^{2}}[\mathbf{d}(u(a))+\mathbf{d}(u(b))].
\end{equation}
\end{lemma}

\begin{proof} Estimate \eqref{1019} is scaling invariant, by rescaling the solution, we only need to prove the estimate with additional 
assumption $\min_{t\in[a,b]}\lambda(t)=1$. In this case, \eqref{1019} can be proved by applying the Fundamental Theorem of Calculus to 
\begin{align}\label{vir_vr}
\partial_{tt}V_R(t)\ge 8 \mathbf{d}(u(t)), \textit{ and } |\partial _t V_R(t)|\lesssim \mathbf{d}(u(t)), \ t\in[a, b]
\end{align}
for some properly chosen $R$. Indeed, the second estimate in \eqref{vir_vr} follows directly from \eqref{gap 1} once $R$ is chosen. To fix this $R$ and control $\partial_{tt} V_R(t)$, we use the fact $\lambda(t)\ge 1$, \eqref{gap 2} and the compactness, in particular \eqref{AA char}, of $u$ to get 
\begin{align*}
|A_R(u(t))|\le8 \mathbf{d}(u(t))
\end{align*}
for some $R=R(\delta_0)$ in both of the two cases $\mathbf{d}(u(t))\ge \delta_0$ and $\mathbf{d}(u(t))< \delta_0$. The estimate on $\partial_{tt}V_R(t)$ follows then quickly from the expression \eqref{1122}. \eqref{1019} is proved. 
\end{proof}

The major obstacle of translating the integration estimate to the point-wise decay of $\dut$ is the uniform lower bound of $\lambda(t)$. We will show this is indeed the case knowing $\dut$ converges to $0$ along a sequence of time, a result that can be deduced again  from Virial analysis.  We prove these results in the following two lemmas.

\begin{lemma}\label{lem:lowerbl}
Let $u$ be the solution of NLS$_{a}$ satisfying \eqref{cond}, there exists a constant $c>0$ such that
\begin{align*}
\inf_{t\in[0,\infty)}\lambda(t)\ge c.
\end{align*}
\end{lemma}
\begin{proof}
Let the sequence $t_n$ be determined by Lemma \ref{L: sub tn} such that $\bd (u(t_n))\to 0$ as $n\to \infty$. There exists sufficiently large $N$ such that 
\begin{align*}
C(\mathbf{d}(u(t_N))+\mathbf{d}(u(t_m)))\le \frac 1{10}, \ \forall m\ge N. 
\end{align*}
where $C$ is the constant in \eqref{intel}. It suffices for us to get the upper bound of $\frac 1{\lambda^2(t)}$ on $[t_N,\infty)$. 

Take any $\tau\in[t_N, \infty)$ and any $m\ge N$ such that $\tau\in [t_N, t_m]$. Applying \eqref{intel} on $[t_N, \tau]$
and Lemma \ref{L: Virial}, we estimate
\begin{align*}
\biggl| \frac 1{\lambda^2(\tau)}-\frac 1{\lambda^2(t_N)}\biggr|&\le C\int_{t_N}^\tau \dut dt\\
&\le C\int_{t_N}^{t_m} \dut dt\\
&\le C\sup_{t\in[t_N, t_m]}\frac 1{\lambda^2(t)}\times \bigl(\mathbf{d}(u(t_N))+\mathbf{d}(u(t_m))\bigr)\\
&\le \frac 1{10} \sup_{t\in[t_N, t_m]}\frac 1{\lambda^2(t)}. 
\end{align*}
Therefore from triangle inequality we have 
\begin{align*}
\frac 1{\lambda^2(\tau)}\le \frac 1{10}\biggl(\sup_{t\in[t_N, t_m]}\frac 1{\lambda^2(t)}\biggr)+\frac 1{\lambda^2(t_N)}, \ \forall \ \tau\in[t_N, \infty). 
\end{align*}
Taking supremum  in $\tau$ on $[t_N, t_m]$ yields
\begin{align*}
\sup_{\tau\in[t_N, t_m]}\frac 1{\lambda^2(\tau)}\le  \frac 2{\lambda^2(t_N)} \text{ and thus } \sup_{\tau\in[t_N, \infty)}\frac 1{\lambda^2(\tau)}\le  \frac 2{\lambda^2(t_N)} 
\end{align*}
by letting $m \to \infty$. The uniform bound for $\tfrac 1{\lambda(\tau)}$ comes from this and the boundedness on the closed interval $[0,t_N]$. Lemma \ref{lem:lowerbl} is proved. 
\end{proof}

Finally we are ready to prove Theorem \ref{thm: subcase}.

\subsection{Proof of Theorem \ref{thm: subcase}}

\begin{proof}

The key of the estimate is to show $\bd(u(t))\to 0$ and in the orthogonal decomposition 
\begin{align}\label{1043} 
\bg_{\theta(t), \mu(t)}u(t)=W+\alpha(t) W+v(t),
\end{align}
all the parameters converges exponentially to their limits. 

We start by considering $\dut$ for which we can use Lemma \ref{L: Virial} and Lemma \ref{lem:lowerbl} to get 
\begin{equation}\label{decay 0}
\int_{t}^{t_{n}}\mathbf{d}(u(s))ds\leq C[\mathbf{d}(u(t))+\mathbf{d}%
(u(t_{n}))].
\end{equation}
Here $\{t_n\}$ is the sequence in Lemma \ref{L: sub tn}, along which $\bd (u(t_n))\to 0$. Taking $t_n\to \infty$ in \eqref{decay 0} gives immediately 
\[\int_t^\infty \bd(u(s)) ds\le C \dut, \ \forall t\ge 0,\]
which together with Gr\"onwall's inequality yields
\begin{equation}
\int_{t}^{\infty }\mathbf{d}(u(s))ds\leq Ce^{-ct},  \label{decay e}
\end{equation}
for some $c, C>0$. 

Now before proving the convergence of $\dut$ we go back to \eqref{intel} and consider the convergence of $\lambda(t)$. Combining the estimates \eqref{intel} and \eqref{decay e} we immediately see that $\frac 1{\lambda^2(t)}$ converges as $t\to \infty$. Therefore, there exists $\lambda_\infty$ such that $\lim_{t\to \infty}\lambda(t)=\lambda_\infty \in (0, \infty]$ from the lower bound estimate Lemma \ref{lem:lowerbl}. To preclude the possibility that $\lambda_\infty=\infty$, we argue by contradiction. Assuming this is the case, i.e.
\begin{equation}\label{820}
\lim_{t\to \infty}\frac 1{\lambda(t)}= 0, 
\end{equation}
and recalling $\bd(u(t_n))\to 0$, for any $\eps>0$, there must exist $N_0\in\mathbb N$ such that 
\begin{equation}\label{1101}
\frac 1{\lambda(t)}<\eps, \ \forall t\ge t_{N_0}\textit{ and }\ \bd(u(t_n))<\eps, \ \forall n\ge N_0. 
\end{equation}
Taking any $t_*\ge t_{N_0}$ and applying \eqref{intel}, \eqref{1019} we obtain
\begin{align*}
\biggl| \frac 1{\lambda^2(t_*)}-\frac 1{\lambda^2(t_n)}\biggr|&\le C\biggl|\int_{t_*}^{t_n} \dut dt\biggr|\le C\int_{t_{N_0}}^{t_n}\dut dt\\
&\le C\sup_{t\in[t_{N_0}, t_n]}\frac 1{\lambda^2(t)}[\bd(u(t_n))+\bd(u(t_{N_0}))]. 
\end{align*}
Letting $n\to \infty$ we have 
\[
\frac 1{\lambda^2(t_*)}\le C\sup_{t\in [t_{N_0},\infty)}\frac 1{\lambda^2(t)} \bd(u(t_{N_0}))\le C\eps\sup_{t\in[t_{N_0},\infty)}\frac 1{\lambda^2(t)}. \]
Choosing $C\eps\le \frac 12$ and taking supremum in $t_*$ over $[t_{N_0},\infty)$, we obtain $\frac 1{\lambda(t)}=0$ for all $t\ge t_{N_0}$, which is a contradiction. Therefore we conclude that 
\begin{equation}\label{841}
\lim_{t\to \infty}\lambda(t)=\lambda_{\infty}\in (0,\infty). 
\end{equation}

Next we turn to proving the convergence 
\begin{align}\label{1115}
\lim_{t\to \infty}\dut=0.
\end{align}
Again we argue by contradiction. If this is not true, there must exist a subsequence in $n$ (which we still use the same notation) and a constant $c\in (0,\delta_0)$ such that 
$\max_{[t_n,t_{n+1}]}\dut\ge c.$ Therefore we can find $\tau_n\in(t_n, t_{n+1})$ such that 
\begin{equation}\label{840}
\bd (u(\tau_n))=c \textit{ and } \dut\le c, \ \forall t\in [t_n, \tau_n]. 
\end{equation}
Applying Lemma \ref{L: Modulation}, integrating $\alpha'$ over $[t_n, \tau_n]$ and using the fact that $\mu=\lambda$, \eqref{841} and \eqref{decay e} we have
\begin{align}
|\alpha(t_n)-\alpha(\tau_n)|&\le \int_{t_n}^{\tau_n}|\alpha'(s)| ds\le C\int_{t_n}^{\tau_n}\frac{|\alpha'(s)|}{\mu^2{(s)}}ds\label{843}\\
&\le C\int_{t_n}^{\tau_n} \bd(u(s)) ds\le Ce^{-ct_n}. \notag
\end{align}
As from Lemma \ref{L: Modulation}, $\alpha(t)\sim \dut$ for $t=t_n,\tau_n$.  Taking $n\to \infty$ in \eqref{843} gives $\alpha(\tau_n)\to 0$, which contradicts with \eqref{840}. Therefore \eqref{1115} is proved. 

Due to \eqref{1115}, orthogonal decomposition remains valid for all large enough $t\ge T_0$. In particular, this implies $\mu(t)=\lambda(t)\to \lambda_\infty$ and $\alpha(t)\sim \dut\to 0$. Combining these estimate and repeating the same estimate in \eqref{843} over the interval $[t,\tau]$ we have 
\[|\alpha(t)-\alpha(\tau)|\le Ce^{-ct},\ \forall t\ge T_0. \]
which by taking $\tau\to \infty$ gives rise to 
\[|\alpha(t)|\le Ce^{-ct}, \ \forall t\ge T_0.\]
From here we apply Lemma \ref{L: Modulation} again to get 
\begin{align}\label{943}
\|v(t)\|_{\dha}+\dut\le C e^{-cT}, \forall t\ge T_0. 
\end{align}
Finally, from the derivatives estimate of $\mu'$, $\theta'$ in Lemma \ref{L: Modulation}, using the boundedness of $\mu(t)$ and \eqref{943}, we know that there exists $\theta_\infty\in  \mathbb{S}^1$ such that
\[|\theta(t)-\theta_\infty|+|\mu(t)-\lambda_\infty|\le Ce^{-ct}, \ \forall t\ge T_0.\]
Therefore finally, we have 
\begin{align*}
\|\bg_{\theta_\infty, \lambda_\infty} u(t)-W\|_{\dha}&=\|u(t)-\bg^{-1} _{\theta_{\infty},\lambda_\infty}W\|_{\dha}\\
&\le \|\bg_{\theta(t),\mu(t)}u(t)-W\|_{\dha}+\|(\bg^{-1}_{\theta(t),\mu(t)}-\bg^{-1}_{\theta_\infty, \lambda_\infty})W\|_{\dha}\\
&\le \alpha(t)\|W\|_{\dha}+\|v(t)\|_{\dha}+C(|\theta(t)-\theta_\infty|+|\mu(t)-\lambda_\infty|)\\
&\le Ce^{-ct}, \ \forall t\ge T_0,
\end{align*}
which by incorporating the finite bound on closed interval $[0,T_0]$ and changing the notation give rise to \eqref{448} in Theorem \ref{thm: subcase}. 

To complete the proof of Theorem \ref{thm: subcase}, we shall prove $\| W_-\|_{\dot S^1 ((-\infty, 0])} <\infty$ by contradiction. Assume $\| W_-\|_{\dot S^1 ((-\infty, 0])} =\infty$ all the above results apply to $W_-(t, x)$ for $t \in (-\infty, 0]$. In particular, Lemma \ref{L: Virial} and Lemma \ref{lem:lowerbl} imply 
\[
\int_{a}^b \bd (W_-(t)) dt \le C \big( \bd(W_-(a)) + \bd(W_-(b))\big), \quad \forall a, b \in \R.
\]
Since $\lim_{|t| \to \infty} \bd(W_-(t)) =0$, we obtain $\int_{-\infty}^\infty \bd (W_-(t)) dt =0$, which implies $W_-\equiv W$. It is a contradiction and the proof is complete. 
\end{proof}

\begin{remark}
1. The above argument verifies b.1) in the statement of Theorem \ref{thm:main2}. By time reversal symmetry, we immediately have b.2). The only missing piece is to show all the solutions on the energy surface with less kinetic energy than that of the ground state must be global solution. There a contradiction argument together with the uniform control on the kinetic energy \eqref{12166} and the compactness of the solution leads to the conclusion, see \cite{R: Kenig focusing} or \cite{R: D Merle} for details. 

2. Except for the compact assumption, the analysis in this chapter is not dimension sensitive and can be extended easily to all  dimensions $d\ge 4$. 
\end{remark}

\section{Exponential convergence in the super-critical case}

In this section, we characterize solutions of $NLS_a$ on the energy surface of $E(W)$ if the kinetic energy is greater than that of the ground state. Being different from Section 7, such solutions do not automatically obey the compactness. We thus add additional  spatial decay and symmetry requirement to get a proper control on the solution. Our result is the following 

\begin{theorem}
\label{thm: super case} 
Let $u$ be a solution to $NLS_a$ 
%defined
satisfying 
\begin{equation}\label{549}
E(u)=E(W),\text{ }\Vert u\Vert _{\dot{H}_{a}^{1}}>\Vert W\Vert _{\dot{H}%
_{a}^{1}},\mbox{ and }u\in H^1_{rad}(\R^3),
\end{equation}
then  the maximal lifespan of $u$ must be finite.
\end{theorem}

%This theorem is proved by first showing that if such $u(t)$ exists  on $[0,\infty )$, then there must be some $(\theta,\mu)\in S^1\times\R^+$ such that 
%\begin{equation} \label{E:unstable-1}
%u(t,x)=e^{-i\theta}\mu^{\frac 12}W_+(\mu^2 t+T, \mu x). 
%\end{equation}
%This together with the fact that $W_+\notin L^2(\R^3)$ in three dimensions from Corollary \ref{cor:411}  immediately {\color{red} proves Theorem \ref{thm: super case}.}

%\begin{theorem}\label{thm:charsuper}
%Let $u$ be a solution of $NLS_a$ satisfying \eqref{549}. Then the maximal life span of $u$ must be finite. 
%\end{theorem}

We start by pointing out some of the implications from the symmetry and regularity assumptions.

\begin{lemma}
\label{lem:mulbsup} Suppose $u$ is the solution in Theorem \ref{thm: super case}. Then we have the following 

1) On the interval $I$ where $\dut<\delta_0$,  there exists $c>0$ such that $\mu(t)$ appearing in the orthogonal decomposition given in Lemma \ref{lem:1st}
\begin{align}\label{604}
\bg_\smt u(t)=W+v(t)
\end{align}
satisfies 
\begin{equation}\label{608}
\mu(t)\ge c, \ \forall t\in I. 
\end{equation}

2) There exists $R_0=R_0(\delta_0, W, \|u\|_2)$ such that when $R\ge R_0$, 
\begin{equation}\label{441}
A_R(u(t))\le \dut. 
\end{equation}

\end{lemma}

\begin{proof}
We first prove \eqref{608}. Taking $L^2$ norm on both sides of \eqref{604} and using $\|v(t)\|_6\lsm \|v\|_{\dha}\le C\delta_0$ from Lemma \ref{L: Modulation} we have
\begin{align*}
\mu (t)\Vert u(t)\Vert _{2}& \geq \Vert W+v(t)\Vert _{L^{2}(|x|\leq 1)} \\
& \geq \Vert W\Vert _{L^{2}(|x|\leq 1)}- C\Vert v(t)\Vert _{6}\geq
\Vert W\Vert _{L^{2}(|x|\leq 1)}-C\delta _{0}.
\end{align*}%
Inequality \eqref{608} then follows from the mass conservation. It is worthwhile to note that in this step that we do not need the radial symmetry. 

We turn to proving \eqref{441}. We first recall the decay estimate for the radial function in three dimensions:
\begin{align*}
|x|^2|u(x)|^2\lesssim \|u\|_2\|\nabla u\|_2,
\end{align*}
which can be proved by using the Fundamental Theorem of Calculus and Hardy's inequality. Inserting this decay estimate into the interpolation, we have
\begin{align}\label{1109}
\|u\|_{L^{6}(|x|\geq R)}&\le \|u\|_2^{\frac 13}\|u\|_{L^{\infty}(|x|\ge R)}^{\frac 23}
\le  R^{-\frac 23}\|u\|_2^{\frac 23}\|\nabla u\|_2^{\frac 13}. 
\end{align}
This estimate together with the first bound in \eqref{gap 3} gives 
\begin{align*}
A_R(u(t))&\lsm \int_{|x|>R}|u(t,x)|^6 dx+\int_{|x|>R}\frac{|u(t,x)|^2}{|x|^2} dx\\
&\lsm R^{-4}\|u\|_2^4(\dut+\|W\|^2_{\dot H^1})+R^{-2}\|u\|_2^2. 
\end{align*}
By taking $R$ large enough depending on $\|u\|_2$, $\delta_0$ and $W$, we immediately have \eqref{441} in the case of $\dut\ge\delta_0$. In the remaining case when $\dut<\delta_0$, \eqref{441} follows directly from \eqref{608} and the second estimate of \eqref{gap 2}. The lemma is proved. 
\end{proof}

We are ready to prove Theorem \ref{thm: super case}. 

\begin{proof}
We argue by contradiction. 
Assume $u(t)$ exists for all $t \in [0, \infty)$, our goal is to show there must be some $(\theta,\mu)\in \mathbb S^1\times\R^+$ such that 
\begin{equation} \label{E:unstable-1}
u(t,x)=e^{-i\theta}\mu^{\frac 12}W_+(\mu^2 t+T, \mu x). 
\end{equation}
As will be explained later this together with the fact that $W_+\notin L^2(\R^3)$ in three dimensions from Corollary \ref{cor:411}  immediately yields a contradiction.

Like in Section \ref{decay_sub}, the key in proving \eqref{E:unstable-1} is to show $\dut\to 0$ and in the orthogonal decomposition 
\begin{equation}\label{154}
\bg_{\smt} u(t)=W+\alpha(t) W+\tilde u(t):=W+v(t),
\end{equation}
all the parameters converge exponentially to their limits. 

We first establish the integral estimate for $\dut$, which again will follow from the Virial analysis.  For $R\ge R_0$, we apply  \eqref{441} to get  
\begin{equation}\label{1144}
\partial _{tt}V_{R}(t)=-16\dut+A_{R}(u(t))\le -15\dut,
\end{equation}%
hence $\partial_t V_R(t)$ decreases on $[0,\infty)$. This further implies that  
\begin{equation}
\partial _{t}V_{R}(t)>0,\ \forall t\geq 0.  \label{Vr mono}
\end{equation}%
Indeed, if this is not true, as $\partial _{t}V_{R}(t)$ is decreasing, there
must exist $t_{0}>0$ such that 
\begin{equation*}
\partial _{t}V_{R}(t)<\partial _{t}V_{R}(t_{0})<0,\ \forall t>t_{0},
\end{equation*}%
which obviously contradicts with the uniform bound $V_R(t) \ge 0$.
%\lsm R^2 \|u\|_2. 

Using the positivity of $\partial_t V_R(t)$ together with the estimate of it from \eqref{gap 1}, we integrate \eqref{1144} over $[t,T]$ to get 
\begin{align*}
\int_t^ T \bd(u(s)) ds\le \partial_t V_R(t)-\partial_t V_R(T)\le \partial_t V_R(t)\le C \dut, \ \ \forall t\ge 0. 
\end{align*}
Taking $T\to \infty$ we obtain 
\begin{equation}\label{1150}
\int_t^{\infty}\bd(u(s)) ds\le C\dut, \; \forall t\ge 0. 
\end{equation}
As a direct implication, there exists a sequence $\{t_n\}\subset(0,\infty)$ such that $\lim_{n\to\infty}\bd(u(t_n))=0$. Therefore we can perform the decomposition \eqref{154} in the neighborhood of $t_n$ for large $n$. We claim that 
\begin{align}\label{1155}
\mu(t_n)\lsm 1. 
\end{align}
Indeed, if this is not true, passing to a subsequence, we have $\mu(t_{n_k})\to \infty$. Along this subsequence we use H\"older and \eqref{154} to estimate 
\begin{align*}
V_{R}(t_{n_k})& =\int_{|x|\leq \eps }\phi
_{R}(x)|u(t_{n_k})|^{2}dx+\int_{|x|>\epsilon }\phi _{R}(x)|u(t_{n_k})|^{2}dx \\
& \leq \eps ^{2}\Vert u(t_{n_k})\Vert _{2}^{2}+R^{4}\Vert u(t_{n_k})\Vert
_{L^{6}(|x|>\eps )}^{\frac13} \\
& \lesssim \eps ^{2}+R^{4}\Vert \bg_{[\theta (t_{n_k}),\mu
(t_{n_k})]}u(t_{n_k})\Vert _{L^{6}(|x|\geq \eps \mu (t_{n_k}))}^{\frac 13}
\\
& \lesssim \eps ^{2}+R^{4}\Vert W\Vert _{L^6(|x|\geq \eps
\mu (t_{n_k}))}^{\frac13}+\Vert v(t_{n_k})\Vert _{6}^{\frac 13}.
\end{align*}%
Taking $n_k\rightarrow \infty $ then $\eps\to 0$,  we obtain $\lim_{n_k\rightarrow
\infty }V_{R}(t_{n_k})=0$ which contradicts \eqref{Vr mono}.

Next, we prove that 
\begin{equation}
\lim_{t\rightarrow \infty }\mathbf{d}(u(t))=0\text{.}  \label{1251}
\end{equation}
We argue by contradiction. If this is not true, there must exist $c\in (0,\delta_0)$, a subsequence in $\{t_n\}$ (for which we use the same notation) and another sequence $\tau_n$ such that
\begin{align}\label{1250}
\tau_n\in(t_n,t_{n+1}), \; \bd(u(\tau_n))=c, \; \dut\in(0, c], \; \forall t\in [t_n,\tau_n].
\end{align}
Take any $t\in[t_n,\tau_n]$, we use the derivative estimate from Lemma \ref{L: Modulation} and \eqref{1150} to obtain
\begin{equation}\label{1015}
\left\vert \frac{1}{\mu (t_{n})^{2}}-\frac{1}{\mu (t)^{2}}\right\vert \lsm
\int_{t_{n}}^t\left\vert \frac{\mu ^{\prime }(t)}{\mu (t)^{3}%
}\right\vert dt\lesssim \int_{t_{n}}^{\infty }\mathbf{d}(u(t))dt\rightarrow
0.
\end{equation}%
This together with the control from \eqref{608} and \eqref{1155} implies 
\begin{equation}\label{1253}
\mu(t)\sim 1, \; \forall t\in [t_n,\tau_n]. 
\end{equation}
Inserting this to the estimate of $\alpha(t)$ we have
\begin{equation} \label{E:alpha-1}
|\alpha(t_n)-\alpha(\tau_n)|\le \int_{t_n}^{\tau_n}|\alpha'(t)|dt\lsm \int_{t_n}^{\tau_n}\frac{|\alpha'(t)|}{\mu^2(t)}\lsm 
\int_{t_n}^\infty \dut \to 0
\end{equation}
as $n\to \infty$. We get a contradiction as $\alpha(t_n)\sim \bd(u(t_n))\to 0$, but $\alpha(\tau_n)\sim \bd(u(\tau_n))\sim 1$ from \eqref{1250}. The convergence of $\dut$ in \eqref{1251} is proved. 

Given \eqref{1251}, we can perform the decomposition for all $t\ge T_0$ and repeat the same argument as in \eqref{1155} to show $\mu(t)\sim 1$. The exponential  convergence of all the parameters follows from the same argument in Section \ref{decay_sub}. We will not repeat here. Therefore \eqref{E:unstable-1} follows from Corollary \ref{cor:411}. However, $W\notin L^2(\R^3)$ together with $W_\pm -W \in L^2(\R^3)$(by \eqref{438})  contradicts with $u(t) \in H^1_{rad}(\R^3)$ and thus Theorem \ref{thm: super case} is proved. 
\end{proof}

\begin{remark}
 Theorem \ref{thm: super case} verifies Theorem \ref{thm:main2} c). As also seen from the proof, the statement in dimension four is the same after a notational change. In dimension five, due to the fact $W^+\in L^2$, any solution obeying \eqref{549} conforms into one of the three scenarios: blowing up both forward and backward in time; coinciding with $W^+$ up to symmetries or $\bar W^+$ up to symmetries. This justifies the remark after Theorem \ref{thm:main2}. 
 \end{remark}

\section{Appendix}

\begin{lemma} [Asymptotic behavior of $G(r)$]\label{lem:asymG}
Let $W$ be the ground state in \eqref{W} and $G(x) = G(|x|) \in \dot H_{rad}^1(\mathbb R^3)$ solving 
\begin{align}\label{eqnG}
G''+\frac 2r G'-\frac{a+2}{r^2}G+5W^4 G=0. 
\end{align}
Then 
\begin{align}\label{asyG}
\mbox{As } r\to 0^+, \ G(r)=O(r^{-\frac 12+\frac 12\sqrt{9+4a}}),\\
\mbox{As } r\to \infty, \ G(r)=O(r^{-\frac 12-\frac 12\sqrt{9+4a}}). 
\end{align}
\end{lemma}
\begin{proof}
We prove the two asymptotics separately. Near $0$, we introduce the new variable 
$$
s=r^\beta,  \ \beta=\sqrt{1+4a},
$$
and rewrite the equation \eqref{eqnG} into 
\begin{align}\label{eqnGs}
G_{ss}+\frac{\beta+1}{\beta s}G_s-\frac{a+2}{\beta^2 s^2}G+\frac {15}{(1+s^2)^{2}} G=0. 
\end{align}
It is easy to see $0$ is the regular-singular point for this ODE with analytic coefficients, therefore there must exist two linear independent solutions in the form of power series: 
\begin{align*}
G_+(s)=s^{\alpha_+}\sum_{n=0}^\infty a_n s^n, \quad a_0=1;\\
G_-(s)=s^{\alpha_-}\sum_{n=0}^\infty b_n s^n, \quad b_0=1. 
\end{align*}
Here, $s^{\alpha_+}$ and $s^{\alpha_-}$ with
$$
\alpha_{\pm}=\frac 1{\beta}(-\frac 12\pm\frac 12\sqrt{9+4a})
$$
 are solutions to Cauchy-Euler equation 
\begin{align*}
G_{ss}+\frac{\beta+1}{\beta s} G_s-\frac{a+2}{\beta^2s^2}G=0. 
\end{align*}
%Inserting $G_+(s)$ into the equation \eqref{eqnGs}, we have 
%\begin{align}\label{ite}
%\sum_{n=0}^\infty a_n\biggl[(n+\alpha_+)&(n+\alpha_+-1)+\frac{\beta+1}\beta(n+\alpha)-\frac{a+2}{\beta^2}\biggr] s^n\\
%&+\frac{15 s^2}{(1+s^2)^2}\sum_{n=0}^\infty a_n s^n=0.
%\end{align}
%Using the Tylor expansion for the analytic function 
%$$
%\frac{15}{(1+s^2)^2}=\sum_{m=0}^\infty c_m s^m,
%$$
%we can write the last term in \eqref{ite} into 
%$$
%s^2\sum_{m=0}^\infty c_m s^m\cdot \sum_{n=0}^\infty a_n s^n=\sum_{n=2}^\infty d_{n-2}s^n, \ \mbox{ where } d_n=\sum_{0\le j\le n}c_j a_{n-j}.
%$$
%Putting all together we have the iterating formula for the coefficient: 
%\begin{align}\label{iteabc}
%\begin{cases}
%&a_n\biggl[(n+\alpha_+)(n+\alpha_+-1)+\frac{\beta+1}\beta(n+\alpha_+)-\frac{a+2}{\beta^2}\biggr]+d_{n-2}=0,\\
%& \ \forall n\ge 0, d_{-1}=d_{-2}=0, \  d_n=\sum_{0\le j\le n}c_j a_{n-j}.
%\end{cases}
%\end{align}
%In particular, $a_0$ is a free constant due to the choice of $\alpha_+$, $a_1=0$ and all $a_n$ with $n\ge 2$ are determined inductively by the equation \eqref{iteabc}. This clearly proves 
Clearly 
$$
G_\pm(s)=O(s^{\alpha_\pm}) \mbox{ as } s\to 0^+.
$$
%By the same argument we can prove 
%$$
%G_-(s)=O(s^{\alpha_-}) \mbox{ as } s\to 0^-. 
%$$
General solutions to \eqref{eqnGs} are 
\[
c_+ G_+(s) + c_- G_-(s).
\]
Since our $G (x) = G(|x|) \in \dot H_{rad}^1 (\mathbb{R}^3)\subset L^6 (\mathbb{R}^3)$, clearly it must hold $c_-=0$ and thus we obtain the desired asymptotics of $G$ near $0$ after we change the variable back to $r$.

%Returning to the $r$ variable and using the condition $G(r)\in \dot H^1$ to eliminate one possibility, we obtain 
%\begin{align*}
%G(r)=O(r^{-\frac 12+\frac 12\sqrt{9+4a}}), \mbox{ as } r\to 0^+.
%\end{align*}

For the asymptotic behavior near infinity, we can reduce the issue into a similar situation by introducing the change of variable 
$$
s=r^{-\beta}.
$$
Equation \eqref{eqnG} in variable $s$ is 
\begin{align}\label{eqnGsl}
G_{ss}+\frac{\beta-1}{\beta s}G_s-\frac{a+2}{\beta^2 s^2}G+\frac {15}{(1+s^2)^2}G=0. 
\end{align}
From a similar analysis, it has two linear independent solutions of order $O(s^{\frac 1\beta(\frac 12\pm\frac 12\sqrt{9+4a})})$ near $s=0$. 
%we see solution to \eqref{eqnGsl} obeys either $G(s)= or $G(s)=O(s^{-\frac 1\beta(\frac 12-\sqrt{9+4a})}$. 
Going back to $r$ variable and using $G(r)\in \dot H^1(\mathbb{R}^3)$, we are able to select the right asymptotics
$$
G(r)=O(r^{-\frac 12-\frac 12\sqrt{9+4a}})
$$
as $r\to \infty$. The Lemma is proved. 
\end{proof}

\begin{lemma}[$\dha$ linear profile decomposition]\label{lem:lpd}

Let $\{f_n\}$ be a bounded sequence in $\dha(\mathbb R^3)$. After passing to a subsequence, there exist $J^*\in\{0,1,2, \cdots\}\cup\{\infty\}$, $\{\phi^j\}_{j=1}^{J^*}\subset \dha(\mathbb R^3)$, $\{(\lambda_n^j, x_n^j)\}_{j=1}^{J^*}\subset \R^+\times\R^3$ such that for every $0\le J\le J^*$, we have the decomposition 
\begin{align*}
f_n=\sum_{j=1}^J\phi_n^j+r_n^J, \ \phi_n^j=(\lambda_n^j)^{-\frac 12}\phi^j\bigl(\frac{x-x_n^j}{\lambda_n^j}\bigr):=g_n^j\phi^j, \ r_n^J\in \dha(\mathbb{R}^3)
\end{align*}
satisfying 
\begin{align}
&\lim_{J\to J^*}\limsup_{n\to \infty}\|r_n^J\|_6=0;\notag\\
&\lim_{n\to \infty}\biggl(\|f_n\|_{\dha}^2-\sum_{j=1}^J\|\phi_n^j\|_{\dha}^2-\|r_n^J\|_{\dha}^2\biggr)=0, \ \forall J;\label{232}\\
&\lim_{n\to \infty}\biggl(\|f_n\|_{6}^6-\sum_{j=1}^J\|\phi_n^j\|_{6}^6-\|r_n^J\|_{6}^6\biggr)=0, \ \forall J.\notag
\end{align}
Moreover, for all $j\neq k$, we have the asymptotic orthogonality property 
\begin{align*}
\lim_{n\to \infty} \biggl(\biggl|\frac{\lambda_n^j}{\lambda_n^k}\biggr|+\biggl|\frac{\lambda_n^k}{\lambda_n^j}\biggr|+\frac{|x_n^j-x_n^k|^2}{\lambda_n^j\lambda_n^k}\biggr)=0. 
\end{align*}
Finally we may also assume for each $j$, either $|x_n^j|/\lambda_n^j\to \infty$ or $x_n^j\equiv0$, therefore 
\begin{align}
\|\phi_n^j\|_{\dha}\to \|\phi^j\|_{X^j}=\begin{cases} 
\|\phi^j\|_{\dot H^1}&\mbox{ as } \frac{|x_n^j|}{\lambda_n^j}\to \infty\\
\|\phi^j\|_{\dha}&\mbox{ as } x_n^j\equiv0. 
\end{cases}
\label{233}
\end{align}

\end{lemma}
\begin{proof}
We use a classical $\dot H^1$ linear profile decomposition developed in the work of G\'{e}rard in \cite{Gerard LPD} as a blackbox to prove this Lemma. For a slight different form we will be using, we refer the readers to see \cite{Clay}. As $\{f_n\}$ is also a bounded sequence in $\dot H^1(\R^3)$, from \cite{Gerard LPD} we obtain a decomposition which enjoys all the properties in Lemma \ref{lem:lpd} except \eqref{232} and \eqref{233}. Convergence \eqref{233} is a result quoted directly from Lemma 3.3 in \cite{R: Monica inverse energy}. Therefore the proof of Lemma \ref{lem:lpd} is reduced to only proving the $\dha$ decoupling \eqref{232} by using all the other statements in this Lemma. Before proving \eqref{232}, we record two properties also coming from the classical result. The first one is what appears in \cite{Gerard LPD} in the position of \eqref{232}, the decoupling in $\dot H^1(\mathbb{R}^3)$:
\begin{align}\label{248}
\lim_{n\to \infty}\biggl(\|f_n\|_{\dot H^1}^2-\sum_{j=1}^J \|\phi^j\|_{\dot H^1}^2-\|r_n^J\|_{\dot H^1}^2\biggr)=0. 
\end{align}
The second one is the weak convergence 
\begin{align}\label{251}
(g_n^j)^{-1}r_n^J\rightharpoonup 0, \;\textit{ weakly in }\dot H^1(\mathbb{R}^3), \;\forall 1\le j\le J. 
\end{align}
In view of \eqref{248} and the expression of $\dha$-norm, we further reduce the matter to proving 
\begin{align}\label{302}
\lim_{n\to \infty}\int_{\R^3}\frac 1{|x|^2}\biggl(|f_n|^2-\sum_{j=1}^J|\phi_n^j|^2-|r_n^J|^2\biggr) dx=0. 
\end{align}
To see \eqref{302}, we use the decomposition to write 
\[|f_n|^2-\sum_{j=1}^J |\phi_n^j|^2-|r_n^J|^2=\sum_{j\neq k}\phi_n^j\bar\phi_n^k+2\Re \sum_{j=1}^J r_n^J\bar \phi_n^j,\]
and estimate the contribution to \eqref{302} from each above term. To estimate the cross term, we write 
\begin{align*}
\int_{\R^3}\frac{\phi_n^j(x)\bar\phi_n^k(x)}{|x|^2}dx&=\int_{\R^3} \phi^j(y)\frac{(g_n^j)^{-1}g_n^k\bar\phi^k(y)}{|y+x_n^j/\lambda_n^j|^2} dy:=A
\end{align*}
and discuss the convergence in two cases. Note here by density argument, we may assume $\phi^j$, $\phi^k\in C_c^\infty(\R^3)$.

 In the first case where $\bigl|\log\frac{\lambda_n^j}{\lambda_n^k}\bigr|\to \infty$, we use Hardy's inequality to obtain 
\begin{align*}
|A&|\le\min\biggl( \biggl\|\frac{\phi^j}{|y+x_n^j/\lnj|^{\frac 34}}\biggr\|_2\biggl\|\frac{(g_n^j)^{-1}g_n^k \phi^k}{|y+x_n^j/\lnj|^{\frac 54}}\biggr\|_2,\; \biggl\|\frac{\phi^j}{|y+x_n^j/\lnj|^{\frac 54}}\biggr\|_2\biggl\|\frac{(g_n^j)^{-1}g_n^k \phi^k}{|y+x_n^j/\lnj|^{\frac 34}}\biggr\|_2\biggr)\\
&\lsm \min\biggl((\lnj/\lnk)^{\frac 14}\|\phi^j\|_{\dot H^{\frac 34}}\|\phi^k\|_{\dot H^{\frac 54}}, \;(\lnj/\lnk)^{-\frac 14}\|\phi^j\|_{\dot H^{\frac 54}}\|\phi^k\|_{\dot H^{\frac 34}}\\
&\lsm\min\bigl( (\lnj/\lnk)^{\frac 14},\;(\lnj/\lnk)^{-\frac 14}\bigr)\to 0,\;\textit{as } n\to \infty. 
\end{align*}
In the second case where $\lnj\sim\lnk$, the orthogonality condition guarantees $\frac{|x_n^j-x_n^k|^2}{\lnj\lnk}\to \infty$ as $n\to \infty$. Going back to the expression $A$, this means the support of $\phi^j$ and $(g_n^j)^{-1}g_n^k\bar\phi^k$ do not overlap, hence $A=0$ for sufficiently large $n$. 

We turn to estimating $\int_{\R^3}\frac{r_n^J(x)\bar \phi_n^j(x)}{|x|^2} dx$, which by changing of variables, can be written as
\begin{align}\label{345}
\int_{\R^3}\frac{r_n^J(x)\bar \phi_n^j(x)}{|x|^2} dx=\int_{\R^3}\frac{(g_n^j)^{-1} r_n^J(y)}{|y+\xnj/\lnj|^2}\bar\phi^j(y) dy.
\end{align}
By density argument we may assume $\phi^j\in C_c^\infty(\R^3/\{0\})$. Recall as part of the classical result, for each $j$ either $\xnj\equiv0$ or $\frac{|x_n^j|}{\lambda_n^j}\to \infty$. In the first case, we immediately have $\eqref{345}\to 0$ from the weak convergence \eqref{251} and the property of $\phi^j$. In the case when $\frac{|x_n^j|}{\lambda_n^j}\to \infty$, assuming $\supp (\phi^j)\subset \{|x|\le R\}$, we can estimate
\begin{align*}
\eqref{345}\le |\xnj/\lnj-R|^{-2}\|g_n^j r_n^J\|_6\|\phi^j\|_{\frac 65}\to 0,
\end{align*}
as $n\to \infty$. Combining all the pieces together we prove \eqref{302}, hence end the proof of Lemma \ref{lem:lpd}. 
\end{proof}

The following lemma is concerned with the modification of the scaling size function $\lambda(t)$ in Section \ref{decay_sub}. Let $u(t)$ be a solution to NLS$_a$ and $\lambda \in C^0([0, \infty), \R^+)$ satisfying \eqref{cond} and \eqref{precomp}.

\begin{lemma} \label{L:lambda-1} (Modification of $\lambda(t)$ in Section \ref{decay_sub}.)
There exist $0< C_1 <C_2$ and a function $\tilde \lambda\in C^0 \big([0, \infty), \R^+\big)$ such that it satisfies \eqref{intel} and $\tilde \lambda'$ exists almost everywhere and 
\begin{equation} \label{E:tlambda}
\frac {\lambda(t)} {\tilde \lambda(t)} \in (C_1, C_2), \quad \forall t\in [0, \infty)
\end{equation}
%satisfies  \in L^\infty $ is  
%Integrating \eqref{perpl} on $[a,b]$ yields
\end{lemma} 

The above property \eqref{E:tlambda} means that the new $\tilde \lambda$ also satisfies \eqref{precomp}.

\begin{proof}  
The proof of this lemma is pure technicality and we divide it into several steps. Let $\delta_0$ be the constant given in Lemma \ref{lem:1st}. 

{\it Step 1.} Let 
\[
A_l = \{t\in (0, \infty) \mid \bd(u(t)) > \frac 23\delta_0\}, \quad A_s = \{t\in (0, \infty) \mid \bd(u(t))\le \frac 23\delta_0\}.
\]
Since $A_l$ is open, it must be the disjoint union of at most countably many intervals 
\[
A_l = \cup_n I_n, \quad I_n =(a_n, b_n), \; a_n < b_n <\infty, \; I_n \cap I_m =\emptyset, \, \forall m\ne n,  
\]
where all $b_n<\infty$ is due to Lemma \ref{L: sub tn}. Recall $\mu(t)$, $t\in A_s$, is the function given in Lemma \ref{lem:1st}. According to Lemma \ref{cop}, 
\[
\frac {\lambda(a_n)}{\mu(a_n)}, \; \frac {\mu(a_n)}{\lambda(a_n)}, \; \frac {\lambda(b_n)}{\mu(b_n)}, \frac {\mu(b_n)}{\lambda(b_n)}, 
\]
are bounded uniformly in $n$. Therefore there exist linear functions $l_n(t)$, $t\in I_n$, bounded uniformly in $n$ such that 
\[
\lambda_1 (t)= \begin{cases} l_n (t) \lambda(t), \qquad & t \in I_n \\ 
\mu (t), & t\in A_s \end{cases}
\]
is continuous. Apparently  $\frac {\lambda_1(t)}{\lambda(t)}$ has positive upper and lower bounds. Therefore $\lambda_1(t)$ satisfies \eqref{precomp} and all the subsequent properties. Moreover $\lambda_1(t)$ is $C^1$ in the interior of $A_s$ and in particular satisfies 
\[
|\lambda_1'(t)| \le C \lambda_1(t)^3 \bd (u(t)), \quad \text{ if } \bd(u(t)) < \frac 23 \delta_0. 
\]

{\it Step 2.} Similarly, let 
\[
B_l = \{t\in (0, \infty) \mid \bd(u(t)) > \frac {\delta_0}3\}, \quad B_s = \{t\in (0, \infty) \mid \bd(u(t))\le \frac {\delta_0}3\}.
\]
Since $B_l$ is open, it must be the disjoint union of at most countably many intervals 
\[
B_l = \cup_n J_n, \quad J_n =(a_n', b_n'), \; a_n' < b_n' <\infty, \; J_n \cap J_m =\emptyset, \, \forall m\ne n. 
\]
We classify the intervals in $B_l$ into two categories by singling out 
\[
\Lambda_l = \{ n \mid \exists \ t \in J_n, \text{ s.t. } \bd(u(t)) \ge\frac 23 \delta_0\}.
\]
We shall only modify $\lambda_1$ in such intervals. 

For any $n \in \Lambda_l$, let 
\[
t_{n*} = \inf \{ t\in J_n \mid \bd(u(t)) \ge \frac 23 \delta_0\}> a_n', \quad t_n^* = \sup  \{ t\in J_n \mid \bd(u(t)) \ge \frac 23 \delta_0\}< b_n'.  
\]
Clearly $\bd(u(t_{n*})) = \bd(u(t_n^*)) = \frac 23\delta_0$. Define $t_0=t_{n*}$ and 
\[
t_{j+1} = t_j + \frac {\eps \delta}{\lambda_1(t_j)^2}, \; j=0, 1, 2 \ldots, \quad k_n = \min\{j \mid t_j >t_{n}^*\}, 
\]
where $\delta$ is given by Lemma \ref{lem:almost constancy} for $\lambda_1$ and $\eps$ is chosen from the next claim.

{\bf Claim.} There exists $\eps\in(0,1]$ such that $k_n<\infty$ and $t_{k_n} \in [t_{k_n-1}, b_n')$. 

In fact, $\lambda_1$ is continuous on $\overline{J_n}$ and thus $\lambda_1^{-2}$ has a positive lower bound, so obviously $k_n <\infty$. To see $t_{k_n} < b_n'$, we argue by contradiction. If this is not true, there must exist a sequence $\eps_m \to 0+$ and intervals such that 
\[b_n'\in[t_{k_n-1},t_{k_n-1}+\frac{\eps_m\delta}{\lambda_1^2(t_{k_n-1})}],\]
which clearly implies $|b_n'-t_{n}^*|\le \frac{\eps_m\delta}{\lambda_1^2(t_{k_n-1})}$. This together with \eqref{modulation 1}, \eqref{modulation 2} and Lemma \ref{lem:almost constancy} gives 
\begin{align*}
|\alpha(b_n')-\alpha(t_n^*)|&\le \sup_{t\in[t_n^*,b_n']}|\alpha'(t)||b_n'-t_n^*|\le \sup_{t\in[t_n^*,b_n']}|\alpha'(t)|\frac{\eps_m\delta}{\lambda_1^2(t_{k_n-1})}\\
&\le C\sup_{t\in[t_n^*,b_n']}\frac{|\alpha'(t)|}{\lambda_1^2(t)}\eps_m\delta\le C\sup_{t\in[t_n^*,b_n']} \dut \eps_m\delta\le C\delta_0\eps_m\delta
\end{align*}
where $\lambda_1 =\mu$ for $t \in [t_n^*, b_n']$ was also used. On the other hand, from \eqref{924}, we can estimate 
\begin{align*}
|\bd(u(b_n'))-\bd(u(t_n^*))|\le 2\|W\|_{\dha}^2|\alpha(b_n')-\alpha(t_n^*)|+C\delta_0^2\le C\delta_0\eps_m\delta+C\delta_0^2\le \frac 1{10}\delta_0. 
\end{align*}
This contradicts with the value of $\bd$ on these two points: $\bd(u(t_n^*))=\frac 23\delta_0$ and $\bd(u(b_n'))=\frac 13\delta_0. $ The claim is proved.

%{\bf Claim.} There exists $\eps \in (0, 1]$ independent of $n$ and $j$ such that $k_n <\infty$ and $t_{k_n} \in (t^*, b_n')$. 
%
%In fact, $\lambda_1$ is continuous on $\overline{J_n}$ and thus $\lambda_1^{-2}$ has a positive lower bound, so obviously $k_n <\infty$. 
%%\[\tau = \inf \{ t \in (t_{k_n-1}, t_{k_n}) \mid \bd(u(t')) < \frac 23\delta_0, \forall t'\in (t, t_{k_n}) \} \in [t_{k_n-1}, t_{k_n}). \]
%From our choice of $t_j$ and Lemma \ref{lem:almost constancy} applied to $\lambda_1$, we have that 
%\[
%\lambda_1 (t) \sim \lambda_1(t_{k_n-1}), \quad t\in [t_{k_n-1}, t_{k_n}], 
%\]
%where the constant depends only on $\lambda_1$. For $t \in [t^*, t_{k_n}] \cap [t^*, b_n']$, we have $\lambda_1(t)=\mu(t)$ and thus \eqref{modulation 1} and \eqref{modulation 2} imply
%\[
%|(\log \alpha(t))'| \lesssim \mu(t)^2 = \lambda_1(t)^{2} \lesssim \lambda_1(t_{k_n-1})^{2} = \frac  {\eps \delta}{t_{k_n} - t_{k_{n-1}}}. 
%\]
%{\color{red} This implies for any $t\in [t_{k_{n-1}},t_{k_n}]$, $|\alpha(t)-\alpha(t_{k_{n-1}})|\lsm \eps \delta\delta_0$. This together with the computation in \eqref{924} shows that 
%\[|\dut-\bd(u(t_{k_{n-1}})|\lsm \eps\delta\delta_0+\delta_0^2\]
%Therefore, by choosing $\eps\in (0, 1]$ and $\delta_0$ reasonably small, we obtain that $\bd(u(t)) > \frac 13 \delta_0$ and so $t_{k_n} < b_n'$. The claim is proved. }

We are ready to start the final modification of $\lambda_1(t)$ on $J_n$ with $n \in \Lambda_l$. For any integer $j\in[0, k_{n-1}]$, there exist constant $\sigma_{n,j,1}$ and $\sigma_{n,j,2}$ such that the function defined by   
\[
\psi_{n, j} (t) = (\sigma_{n, j, 1} t + \sigma_{n, j, 2})^{-\frac 12},\]
satisfies
\[ \psi_{n, j} (t_j) = \lambda_1 (t_j), \; \psi_{n, j} (t_{j+1}) = \lambda_1 (t_{j+1}).\]

Since $\lambda_1(t_j) \sim \lambda_1(t_{j+1})$ and $\psi_{n, j}$ is monotonic, its boundary condition implies 
\[
\psi_{n, j} \sim \lambda_1 \; \text{ on } [t_j, t_{j+1}].
\]
One may compute explicitly 
\[
|\psi_{n, j}'(t)|/\psi_{n, j}(t)^3 = \frac 12 |\sigma_{n, j, 1}| = \frac {|\lambda_1(t_{j+1})^{-2} -\lambda_1(t_{j})^{-2}|}{2(t_{j+1} - t_j)} \lesssim \frac {\lambda_1(t_{j})^{-2}}{t_{j+1} - t_j} = (\eps \delta)^{-1}. 
\]
Define 
\[
\tilde \lambda(t) = \begin{cases} \psi_{n, j} (t), \qquad & t \in [t_j, t_{j+1}] \subset J_n \subset B_l, \; n \in \Lambda_l, \; 0\le j < k_n, \\ 
\lambda_1 (t), & \text{otherwise}. \end{cases}
\]
Clear $\tilde \lambda$ satisfies \eqref{intel} as locally it is equal to $\psi_{n, j}$ or $\mu$ both of which satisfy \eqref{intel}. The construction also ensures $\lambda(t) \sim \tilde \lambda(t)$. 
\end{proof}

%\subsection*{Acknowledgement}
%
%K. Yang was partially supported by the Ballard Seashore Dissertation Fellowship from
%University of Iowa.

\end{document}